\documentclass[12pt,a4paper,fleqn]{article}
\usepackage{a4wide,amsfonts,amsmath,latexsym,amssymb,euscript,graphicx,units,mathrsfs}

\usepackage{pstricks, pst-plot, pst-node, pst-tree}
\usepackage{graphicx}
\usepackage{color}
\usepackage{amssymb}
\usepackage{amssymb}
\usepackage[T1]{fontenc}
\usepackage{latexsym}
\usepackage{xypic}
\usepackage{eufrak}
\usepackage{euscript}
\usepackage{amsfonts,amsmath}
\usepackage{verbatim}
\usepackage{fancyhdr}
\usepackage[english]{babel}
\usepackage{mathrsfs}
\usepackage{units}

\newtheorem{prop}{Proposition}[section]
\newtheorem{cor}[prop]{Corollary}
\newtheorem{corollary}[prop]{Corollary}
\newtheorem{lemme}[prop]{Lemma}
\newtheorem{lemma}[prop]{Lemma}
\newtheorem{rem}[prop]{Remark}
\newtheorem{remark}[prop]{Remark}
\newtheorem{thm}[prop]{Theorem}
\newtheorem{theorem}[prop]{Theorem}
\newtheorem{defi}[prop]{Definition}

\renewcommand{\geq}{\geqslant}
\def\leq{\leqslant}

\newcommand{\R}{\mathbb{R}}

\def\1{{\mathbf{1}}}

\def\1{{\mathbf{1}}}
\def\0.5{{\frac{1}{2}}}

\newcommand{\fin}
{ \vspace{-0.6cm}
\begin{flushright}
\mbox{$\Box$}
\end{flushright}
\noindent }

\newcommand{\qed}{\nopagebreak\hspace*{\fill}
{\vrule width6pt height6ptdepth0pt}\par}

\begin{document}

\begin{center}
{\large{\bf Universal Gaussian fluctuations of non-Hermitian matrix ensembles: from weak convergence to almost sure CLTs}}\\~\\
by Ivan Nourdin\footnote{Laboratoire de Probabilit{\'e}s et
Mod{\`e}les Al{\'e}atoires, Universit{\'e} Pierre et Marie Curie,
Bo{\^\i}te courrier 188, 4 Place Jussieu, 75252 Paris Cedex 5,
France. Email: {\tt ivan.nourdin@upmc.fr}} and Giovanni
Peccati\footnote{Unit\'{e} de Recherche en Math\'{e}matiques, Universit\'{e} du Luxembourg,
162A, avenue de la Fa\"{\i}encerie, L-1511 Luxembourg, Grand-Duchy of Luxembourg. Email: \texttt{giovanni.peccati@gmail.com}}\\
{\it Universit\'e Paris VI and Universit\'e du Luxembourg}\\~\\
\end{center}
{\small \noindent {\bf Abstract.}
In the paper \cite{Noupecrei3}, written in collaboration with Gesine Reinert, we proved a universality
principle for the Gaussian Wiener chaos. In the present work,
we aim at providing an original example of application of this principle in the framework
of random matrix theory.
More specifically, by combining the result in \cite{Noupecrei3} with
some combinatorial estimates,
we are able to prove multi-dimensional central limit theorems for the spectral moments (of arbitrary degrees)
associated
with random matrices with real-valued i.i.d. entries, satisfying some appropriate moment conditions.
Our approach has the advantage of yielding, without extra effort,
bounds over classes of smooth (i.e., thrice differentiable) functions, and it
 allows to deal directly with discrete distributions. As a further application of our estimates, we provide a new ``almost sure central limit theorem'', involving logarithmic means of functions of vectors of traces.
\\

\noindent {\bf Key words:} Almost sure central limit theorems; Central limit theorems; Eigenvalues; Fourth-moment criteria;
Invariance principles; Non-Hermitian random matrices; Normal approximation; Spectral moments;
Universality.\\

\noindent
{\bf 2000 Mathematics Subject Classification:} 60F05; 60G15; 60H05; 60H07. }

\tableofcontents

\section{Introduction}

\subsection{Overview and main results}

In the paper \cite{Noupecrei3}, written in collaboration with Gesine Reinert, we proved several
{\sl universality results}, involving sequences of random vectors whose components have the form of finite
homogeneous sums based on sequences of independent random variables. Roughly speaking, our main finding
implied that, in order to study the normal approximations of homogeneous sums (and under
suitable moment conditions) it is always possible to replace the original sequence with an i.i.d. Gaussian
family. The power of this approach resides in the fact that homogeneous sums associated with Gaussian sequences
are indeed elements of the so-called {\sl Wiener chaos}, so that normal approximations can be established by
means of the general techniques developed in \cite{NP-PTRF, nunugio, PTu04} -- that are based on a powerful
interaction between standard Gaussian analysis, {\sl Malliavin calculus} (see e.g. \cite{nualartbook}) and
{\sl Stein's method} (see e.g. \cite{chen-shao}). Moreover,
in the process one always recovers uniform bounds over suitable classes of smooth functions.

The aim of this paper is to introduce these techniques into the realm of random matrix theory. More specifically,
our goal is to use the universality principles developed in \cite{Noupecrei3}, in order to prove
the forthcoming Theorem \ref{Main}, which consists in a multidimensional
central limit theorem (CLT) for traces of non-Hermitian random matrices with i.i.d. real-valued entries. As explained below, the computations and estimates involved in the proof of Theorem \ref{Main} will be further applied in Section \ref{S : ASCLT}, where we will establish an {\sl almost sure central limit theorem} (ASCLT) for logarithmic means associated with powers of large non-Hermitian random matrices. See Theorem \ref{thm-ASCLT} for a precise statement -- as well as \cite{hormann} for a general discussion on ASCLTs.

\medskip

Now let $X$ be a centered real random variable, having unit variance and with finite moments of all orders, that is, $E(X)=0$, $E(X^2)=1$ and $E|X|^n<\infty$ for every $n\geq 3$.
We consider a doubly indexed collection ${\bf X} =\{X_{ij} : i,j\geq 1\}$
of i.i.d. copies of $X$. For every integer $N\geq 2$, we denote by $X_{N}$ the $N\times N$ random matrix
\begin{equation}\label{Def : RM}
X_{N} =\left \{\frac{X_{ij}}{\sqrt{N}} : i,j=1,...,N \right\},
\end{equation}
and by ${\rm Tr}(\cdot)$ and $X^k_N$, respectively, the usual trace operator and the $k$th power of $X_N$.

\begin{theorem}\label{Main}
Let the above notation prevail. Fix $m\geq 1$, as well as
integers \[1\leq k_1 <\ldots<k_m.\]
Then, the following holds.
\begin{itemize}
\item[\rm (i)] As $N\rightarrow\infty$,
\begin{equation}\label{Main CLT}
\Big({\rm Tr}(X^{k_1}_{N})-E\left[{\rm Tr}(X^{k_1}_{N})\right]
,\ldots,
{\rm Tr}(X^{k_m}_{N})
-E\left[{\rm Tr}(X^{k_m}_{N})\right]\Big)
\stackrel {\rm Law}{\longrightarrow}  \big(Z_{k_1} ,\!..., Z_{k_m}\big),
\end{equation}
where ${\bf Z}=\{Z_k : k\geq 1\}$ denotes a collection of real independent centered Gaussian random variables such that, for every $k\geq 1$, $E(Z_k^2) =k$.
\item[\rm (ii)] Write $\beta = E|X|^3$. Suppose that the function $\varphi : \R^m \rightarrow \R $ is
thrice differentiable and that its partial derivatives up to the order three are bounded by some
constant $B<\infty$. Then, there exists a finite
constant $C =C(\beta,B,m,k_1,...,k_m)$, not depending on $N$, such that
    \begin{eqnarray}\label{BerryBerry}
&&    \Bigg|
E\left[\varphi\left(
\frac{{\rm Tr}(X^{k_1}_{N})-E[{\rm Tr}(X^{k_1}_{N})]}{\sqrt{{\rm Var}({\rm Tr}(X^{k_1}_{N}))}}
,\ldots,
\frac{{\rm Tr}(X^{k_m}_{N})-E[{\rm Tr}(X^{k_m}_{N})]}{\sqrt{{\rm Var}({\rm Tr}(X^{k_m}_{N}))}}]
\right)\right]\\
&&\hskip7cm
- E\left[\varphi\left(\frac{Z_{k_1}}{\sqrt{k_1}},...,\frac{Z_{k_m}}{\sqrt{k_m}}\right)\right]\Bigg|
\leq C\,N^{-1/4}.\notag
    \end{eqnarray}
\end{itemize}
\end{theorem}

\begin{remark}\label{rkrkrk}
{\rm
\begin{enumerate}
\item We chose to state and prove Theorem \ref{Main} in the case of
non-Hermitian matrices with {\sl real-valued} entries, mainly in order to facilitate the connection with
the universality results proved in \cite{Noupecrei3}. However, our techniques may be extended
to the case where the random variable $X$ is {\sl complex-valued} and with finite absolute moments of every order.
This line of research will be pursued elsewhere.
One should also note that, differently from \cite{RiderSilvAop2006}, in the present paper we do not use any
technique coming from complex analysis.
\item Fix an integer $K\geq 2$ and assume that $E|X|^{2K}<\infty$,
while higher moments are allowed to be possibly infinite. By inspection of the forthcoming
proof of Theorem \ref{Main}, one sees that the CLT (\ref{Main CLT}) as well as the bound (\ref{BerryBerry})
continue to hold,
as long as the integers $k_1,...,k_m$ verify $k_j \leq K$ for $j=1,...,m$.
\item In a similar vein as at the previous point, by imposing adequate uniform bounds on moments one can
easily adapt our techniques in order to deal with random matrices whose entries are independent but not
identically distributed. One crucial fact supporting this claim is that the universality principles of
Section \ref{S : UNivers} hold for collections of independent, and not necessarily identically distributed,
random variables.
\item For non-Hermitian matrices, limits of moments are not sufficient to
provide an exhaustive description of the limiting
spectral measure or of the fluctuations around it. Rather, one would need to consider polynomials in the
eigenvalues and their complex conjugates. These quantities
cannot be represented using traces of powers of $X_N$, so that our approach cannot be extended to
this case.
\end{enumerate}
}
\end{remark}

\subsection{Discussion}

In this section we compare our Theorem \ref{Main}
with some related results proved in the existing probabilistic literature.\\

{\bf 1.} In the paper \cite{RiderSilvAop2006}, Rider and Silverstein proved
the following CLT.

\begin{thm}\label{ridersilver}
Let $X$ be a {\rm complex}
random variable such that
$E(X)=E(X^2)=0$, $E(|X|^2)=1$, $E(|X|^k)\leq k^{\alpha k}$, $k\geq 3$
(for some $\alpha>0$) and ${\rm Re}(X)$, ${\rm Im}(X)$ possess a joint bounded density.
For $N\geq 2$, let $X_N$ be defined as in (\ref{Def : RM}).
Consider the space $\mathcal{H}$ of functions $f:\mathbb{C}\to\mathbb{C}$ which are analytic
in a neighborhood of the disk $|z|\leq 4$ and otherwise bounded. Then, as $N\to\infty$,
the random field
\[
\{{\rm Tr}(f(X_{N}))-E\left[{\rm Tr}(f(X_{N}))\right]:\,f\in\mathcal{H}\}
\]
converges in the sense of finite-dimensional distributions (f.d.d.) to the centered complex-valued
Gaussian field
$
\{Z(f):\,f\in\mathcal{H}\},
$
whose covariance structure is given by
\begin{equation}\label{consistent}
E[Z(f)\overline{Z(g)}]=\int_{\mathbb{U}}f'(z)\overline{g'(z)}\frac{d^2z}{\pi}.
\end{equation}
Here, $\mathbb{U}=\{z\in \mathbb{C} : |z|\leq 1\}$ is the unit disk,
and $d^2z/\pi$ stands for the uniform measure on $\mathbb{U}$
(in other words, $d^2z=dxdy$ for $x,y\in\R$ such that $z=x+iy$).
\end{thm}
By using the elementary relations: for every integers $n,m\geq 0$,
\begin{equation*}
\frac1\pi\int_{\mathbb{U}} z^n \overline{z}^m d^2z =\left\{
\begin{array}{ll}
(n+1)^{-1} & \text{if } m=n \\
0 & \text{otherwise,}
\end{array}
\right.
\end{equation*}
one sees that our Theorem \ref{Main} can be reformulated by  saying that
\begin{equation}\label{reformulation}
\{{\rm Tr}(f(X_{N}))-E\left[{\rm Tr}(f(X_{N}))\right]:\,f\in{\rm Pol}(\mathbb{C})\}
\overset{\rm f.d.d.}{\longrightarrow}
\{Z(f):\,f\in{\rm Pol}(\mathbb{C})\},
\end{equation}
where the covariance structure of
$\{Z(f):\,f\in{\rm Pol}(\mathbb{C})\}$ is given by (\ref{consistent}).
It follows that Theorem \ref{Main} {\sl roughly} agrees with Theorem \ref{ridersilver}.
however, we stress that the framework
of \cite{RiderSilvAop2006} is different from ours, since the findings therein cannot be applied
to the real case due to the assumption that real and imaginary parts of entries must possess
a joint bounded density. In addition, also note that (differently from \cite{RiderSilvAop2006}) we do not introduce in the present paper
any requirement on the absolute continuity of the law of the real random variable $X$, so that the framework of
our Theorem
\ref{Main} contemplates every discrete random variable with values in a finite set and with
unit variance. \\

{\bf 2.} One should of course compare the results of this paper with the CLTs involving traces of
{\sl Hermitian}
random matrices, like for instance Wigner random matrices. One general reference in this direction
is the fundamental paper by Anderson and Zeitouni \cite{AZ}, where the authors obtain CLTs for traces
associated with large classes of (symmetric) band matrix ensembles, using a version of the classical
method of moments based on graph enumerations. It is plausible that some of the findings of the present
paper could be also deduced from a suitable extension of the combinatorial devices introduced in \cite{AZ}
to the case of non-Hermitian matrices. However, proving Theorem \ref{Main} using this kind of techniques
would require estimates for arbitrary joint moments of traces, whereas our approach merely requires the
computation of variances and fourth moments. Also, the findings of \cite{AZ} do not allow to directly
deduce bounds such as (\ref{BerryBerry}). We refer the reader e.g. to Guionnet
\cite{GuionnetBook} or to Anderson {\it et al.} \cite{AGZ}, and the references therein, for a detailed overview of existing asymptotic results
for large Hermitian random matrices.  \\

{\bf 3.} The general statement proved by Chatterjee in \cite[Theorem 3.1]{Chatterjee_ptrf} concerns the normal
approximation of linear statistics of random matrices that are possibly non-Hermitian. However, the techniques
used by the author require that the entries can be re-written as smooth transformations of Gaussian random
variables. In particular, the findings of \cite{Chatterjee_ptrf} do not apply to discrete distributions.
On the other hand,
the results of \cite{Chatterjee_ptrf} also provide uniform bounds (based on Poincar\'{e}-type
inequalities and in the total variation distance) for one-dimensional CLTs.
Here, we do { not} introduce
any requirements on the absolute continuity of the law of the real random variable $X$,
and we get bounds for {\sl multi}-dimensional CLTs. \\

{\bf 4.} Let us denote by $\{\lambda_j (N) : j=1,...,N\}$ the complex-valued (random)
eigenvalues of $X_{N}$, repeated according to their multiplicities.
Theorem \ref{Main} deals with the spectral moments of $X_{N}$, that are defined by the relations:
\begin{equation}\label{spectralMom}
N\times\int z^k d\mu_{X_{N}}(z) = \sum_{j=1}^N \lambda_j(N)^k = {\rm Tr}(X^k_{N}),\quad N\geq 2,\quad
k\geq 1,
\end{equation}
where $\mu_{X_{N}}$ denote the spectral measure of $X_N$. Recall that
\begin{equation}\label{Eq : SpMeas}
\mu_{X_N}(\cdot) = \frac{1}{N}\sum_{j=1}^N \delta_{\lambda_j(N)}(\cdot),
\end{equation}
where $\delta_{z}(\cdot)$ denotes the Dirac mass at $z$, and observe that one has also the alternate expression
\begin{equation}\label{trace=sum}
{\rm Tr}(X^k_{N})= N^{-\frac{k}{2}} \sum_{i_1,...,i_k=1}^N X_{i_1 i_2}X_{i_2 i_3}\cdot\cdot\cdot X_{i_k i_1}.
\end{equation}
It follows that our Theorem \ref{Main} can be seen as a partial (see Remark \ref{rkrkrk} (4) above)
characterization of the Gaussian fluctuations associated with the so-called {\sl circular law},
whose most general version has been
recently proved by Tao and Vu:

\begin{thm}[Circular law, see \cite{TV}]\label{T : CircLaw} Let $X$ be a complex-valued random variable, with mean zero and unit variance.
For $N\geq 2$, let $X_N$ be defined as in (\ref{Def : RM}). Then, as $N\rightarrow\infty$,
the spectral measure $\mu_{X_N}$ converges almost surely to the uniform measure
on the unit disk $\mathbb{U}=\{z\in \mathbb{C} : |z|\leq 1\}$. The convergence takes place
in the sense of the vague topology.
\end{thm}

To see why Theorem \ref{Main} concerns fluctuations around the circular law, one can proceed as follows. First observe that, since $E(X^2)=1$ and $E(X^4 )<\infty$ by assumption, one can use a result by Bai and Yin \cite[Theorem 2.2]{BaiYin} stating that, with probability one,
    \begin{equation}\label{EQ : baiYin}
    \limsup_{N\rightarrow\infty} \max_{j=1,...,N}|\lambda_j(N)|\leq 1.
    \end{equation}
    Now fix a polynomial $p(z)$. Elementary considerations yield that, since (\ref{EQ : baiYin}) and the circular law are in order, with probability one
    \begin{equation}\label{LGN}
    \frac{1}{N}{\rm Tr}(p(X_{N})) \rightarrow \frac1\pi \int_{\mathbb{U}} p(z) d^2z = p(0).
    \end{equation}
On the other hand, it is not difficult to see that, for every $k\geq 1$ and as $N\to\infty$,
\[
E\left[\int z^k d\mu_{X_N}(z)\right]=E\left[\frac1N {\rm Tr}(X_N^k)\right]\to 0
\]
(one can use e.g. the same arguments exploited in the second part of the proof Proposition \ref{1dim} below).
This implies in particular, for every complex polynomial $p$,
\begin{equation}\label{nomoments2}
E\left[\frac1N {\rm Tr}(p(X_n))\right]\to p(0)=\frac{1}{\pi}\int_{\mathbb{U}}p(z)d^2z.
\end{equation}
    By (\ref{LGN}) and (\ref{nomoments2}), one has therefore that the quantities $\frac{1}{N}{\rm Tr}(p(X_{N}))$
and $E(\frac{1}{N}{\rm Tr}(p(X_{N})))$ both converge to $p(0)$, and
(\ref{reformulation}) ensures that, for $N$ sufficiently large, the difference \[ {\rm Tr}(p(X_{N}))-Np(0)-\left[ E\left({\rm Tr}(p(X_{N}))\right)-Np(0)\right]\] has approximately a centered Gaussian distribution with variance
$\frac1\pi
\int_{\mathbb{U}} |p'(z)|^2 d^2z.$
Equivalently, one can say that the random variable $\frac{1}{N}{\rm Tr}(p(X_{N}))$
tends to concentrate around its mean as $N$ goes to infinity, and (\ref{reformulation})
describes the Gaussian fluctuations associated with this phenomenon.

On the other hand, one crucial feature of the proof of the circular law provided in \cite{TV} is that it is based
on a universality principle.
This result basically states that, under adequate conditions, the distance between the spectral measures of
(possibly perturbed) non-Hermitian matrices converges systematically to zero,
so that Theorem \ref{T : CircLaw} can be established by simply focussing on the case where $X$ is complex
Gaussian (this is the so-called Ginibre matrix ensemble, first introduced in \cite{Ginibre1965}).
It is interesting to note that our proof of Theorem \ref{Main} is also based on a universality result.
Indeed, we shall show that the relevant part of the vector on the LHS of (\ref{Main CLT}) (that is, the part
not vanishing at infinity) has the form of a collection of homogeneous sums with fixed orders.
This implies that the CLT in (\ref{Main CLT}) can be deduced from the results established in \cite{Noupecrei3},
where it is proved that the Gaussian Wiener chaos has a universal character with respect to Gaussian
approximations. Roughly speaking, this means that, in order to prove a CLT for a vector of general homogeneous
sums, it is sufficient to consider the case where the summands are built from an i.i.d. Gaussian sequence.
This phenomenon can be seen as a further instance of the so-called Lindeberg invariance principle
for probabilistic approximations, and stems from powerful approximation results by Rotar' \cite{Rotar2}
and Mossel {\it et al.} \cite{MOO}. See the forthcoming Section \ref{S : UNivers} for precise statements.\\

{\bf 5}.
We finish this section by listing and discussing very briefly some other results related to Theorem \ref{Main},
taken from the existing probabilistic literature.
\begin{itemize}
\item[-] In Rider \cite{riderptrf2004} (but see also Forrester \cite{forrester1999}), one can find a CLT for (possibly discontinuous) linear statistics of the eigenvalues associated with complex random matrices in the Ginibre ensemble. This partially builds on previous findings by Costin and Lebowitz \cite{CostLeb}.
\item[-] Reference \cite{RiderVirag}, by Rider and Virag, provides further insights into limit theorems involving sequences in the complex Ginibre ensemble. In particular, one sees that relaxing the assumption of analyticity on test functions yields a striking decomposition of the variance of the limiting noise, into the sum of a ``bulk'' and of a ``boundary'' term.  Another finding in \cite{RiderVirag} is an asymptotic characterization of characteristic polynomials, in terms of the so-called {\sl Gaussian free field}.
\item[-] Finally, one should note that the Gaussian sequence ${\bf Z}$ in Theorem \ref{Main} also appears
when dealing with Gaussian fluctutations of vectors of traces associated with large, Haar-distributed unitary
random matrices. See e.g. \cite{DE}
and \cite{DS} for two classic references on the subject.
\end{itemize}

\subsection{Proof of Theorem \ref{Main}: the strategy}\label{strategy}
In order to prove (\ref{Main CLT}) (and (\ref{BerryBerry}) as well), we use an original combination of techniques, which are based both on the universality results of \cite{Noupecrei3} and on combinatorial considerations. The aim of this section is to provide a brief outline of this strategy.

For $N\geq 1$, write $[N]= \{1,...,N\}$. For $k\geq 2$, let us denote by $D_N^{(k)}$ the collection
of all vectors ${\bf i}=(i_1,\ldots,i_k)\in[N]^k$ such that all pairs $(i_a,i_{a+1})$, $a=1,\ldots,k$, are different (with the convention that $i_{k+1}=i_1$), that is, ${\bf i}\in D_N^{(k)}$ if and only if $(i_a,i_{a+1})\neq (i_b,i_{b+1})$ for every $a\neq b$.
Now consider the representation given in (\ref{trace=sum}) and, after subtracting the expectation, rewrite the resulting expression as follows:
\begin{eqnarray}\label{representationINTRO}
&&{\rm Tr}(X_N^k) - E\left[
{\rm Tr}(X_N^k)
\right] \notag\\
&=&
N^{-\frac{k}{2}}  \sum_{i_1,...,i_k=1}^N
\big(X_{i_1 i_2}X_{i_2 i_3}\cdot\cdot\cdot X_{i_k i_1}
 - E[
X_{i_1 i_2}X_{i_2 i_3}\cdot\cdot\cdot X_{i_k i_1}
]\big)
\\
&=&
N^{-\frac{k}{2}}  \sum_{{\bf i}\in D_N^{(k)}} X_{i_1 i_2}X_{i_2 i_3}\cdot\cdot\cdot X_{i_k i_1}\notag\\
&&\hskip2cm
+N^{-\frac{k}{2}}  \sum_{{\bf i}\not\in D_N^{(k)}}
\big(X_{i_1 i_2}X_{i_2 i_3}\cdot\cdot\cdot X_{i_k i_1}
 - E[
X_{i_1 i_2}X_{i_2 i_3}\cdot\cdot\cdot X_{i_k i_1}
]\big).
 \label{diagonals}
\end{eqnarray}
Our proof of (\ref{Main CLT})  is based on the representation (\ref{representationINTRO})--(\ref{diagonals}), and it is divided in two (almost independent) parts.\\
\\
{\bf I.} In Section \ref{section3}, we shall prove that the following multi-dimensional CLT takes place for every integers $2\leq k_1<...<k_m$:
\begin{eqnarray}\label{todo}
&&\left(N^{-1/2}\sum_{i=1}^NX_{ii},\,\,\,
N^{-\frac{k_1}{2}}  \sum_{{\bf i}\in D_N^{(k_1)}} X_{i_1 i_2}X_{i_2 i_3}\cdot\cdot\cdot X_{i_{k_1} i_1}
,\,\ldots \right. \\
&&\quad\quad\quad\quad\quad\quad\quad\quad \quad \left. \ldots,\,
N^{-\frac{k_m}{2}}  \sum_{{\bf i}\in D_N^{(k_m)}} X_{i_1 i_2}X_{i_2 i_3}\cdot\cdot\cdot X_{i_{k_m} i_1}
\right)
\!\stackrel{\rm Law}{\longrightarrow} \! \big(Z_1,Z_{k_1} ,..., Z_{k_m}\big),\notag
\end{eqnarray}
for ${\bf Z}=\{Z_i:\,i\geq 1\}$ as in Theorem \ref{Main}.
In order to prove (\ref{todo}), we apply the universality result obtained in \cite{Noupecrei3} (and stated in a convenient form in the subsequent Section \ref{S : UNivers}). This result roughly states that, in order to show (\ref{todo}) in full generality, it is sufficient to consider the special case where the collection ${\bf X} = \{X_{ij} : i,j\geq 1\}$ is replaced by an i.i.d. centered Gaussian family ${\bf G} = \{G_{ij} : i,j\geq 1\}$, whose elements have unit variance. In this way, the components of the vector on the LHS of (\ref{todo}) become elements of the so-called {\sl Gaussian Wiener chaos} associated with ${\bf G}$:
it follows that one can establish the required CLT by using the general criteria for normal approximations on a fixed Wiener chaos, recently proved in \cite{NP-PTRF, nunugio, PTu04}. Note that the results of \cite{NP-PTRF, nunugio, PTu04} can be described as a ``simplified method of moments'': in particular, the proof of (\ref{todo}) will require the mere computation of quantities having the same level of complexity of covariances and fourth moments.\\
\\
{\bf II.} In Section \ref{section4}, we shall prove that the term (\ref{diagonals}) vanishes as $N\rightarrow\infty$, that is, for every $k\geq 2$,
\begin{eqnarray}\notag
R_N(k) :=N^{-\frac{k}{2}}  \sum_{{\bf i}\not\in D_N^{(k)}}
\big(X_{i_1 i_2}X_{i_2 i_3}\cdot\cdot\cdot X_{i_k i_1}
 - E[
X_{i_1 i_2}X_{i_2 i_3}\cdot\cdot\cdot X_{i_k i_1}
]\big)\to 0\quad\mbox{in $L^2(\Omega)$}.\\
\label{remainder}
\end{eqnarray}
The proof of (\ref{remainder}) requires some subtle combinatorial analysis, that we will illustrate by means of graphical devices, known as {\sl diagrams}. Some of the combinatorial arguments and ideas developed in Section \ref{section4} should be compared with the two works by Geman \cite{Geman80, Geman1986}.\\

Then, the upper bound (\ref{BerryBerry}) will be deduced in Section \ref{SS : Last Word} from the estimates obtained at the previous steps. \\

\subsection{An application to almost sure central limit theorems}

As already pointed out, one of the main advantages of our approach is that it yields explicit estimates for the normal approximation of vectors of traces of large random matrices -- see e.g. relation (\ref{BerryBerry}). In Section \ref{S : ASCLT}, we shall show that these estimates can be effectively used in order to deduce multivariate almost sure central limit theorems (ASCLTs), such as the one stated in the forthcoming Theorem \ref{thm-ASCLT}. In particular, this result involves powers of non-Hermitian random matrices and sheds further light on the asymptotic behavior of their traces.
To the best of our knowledge, Theorem \ref{thm-ASCLT} is the first ASCLT ever proved in the
context of traces of random matrices.

\begin{thm}\label{thm-ASCLT}
Fix $m\geq 1$, as well as
integers $k_m>\ldots>k_1\geq 1$, and let the Gaussian vector $(Z_{k_1},\ldots,Z_{k_m})$ be defined as in Theorem \ref{Main}.
Then, a.s.-$P$,
\begin{equation}\label{Main-ASCLT}
\frac{1}{\log N}\sum_{n=1}^N
\frac{1}{n}\varphi\Big({\rm Tr}(A^{k_1}_{n})-E\left[{\rm Tr}(A^{k_1}_{n})\right]
,\ldots,
{\rm Tr}(A^{k_m}_{n})
-E\left[{\rm Tr}(A^{k_m}_{n})\right]\Big)
{\rightarrow}  E\big[\varphi\big(Z_{k_1} ,\!..., Z_{k_m}\big)\big],
\end{equation}
as $N\rightarrow\infty$, for every continuous and bounded function $\varphi:\R^m\to\R$.
\end{thm}

\begin{rem}{\rm
\begin{enumerate}
\item Fix $m\geq 1$ and, for every $N\geq 1$, denote by $\rho_{N}$ the discrete random measure on $\R^m$ assigning mass $(n \log(N))^{-1}$ to the points \[\left({\rm Tr}(A^{k_1}_{n})-E\left[{\rm Tr}(A^{k_1}_{n})\right]
,\ldots,
{\rm Tr}(A^{k_m}_{n})
-E\left[{\rm Tr}(A^{k_m}_{n})\right]\right), \quad  n=1,...,N. \]
Then, the usual characterization of weak convergence imply that relation (\ref{Main-ASCLT}) is indeed equivalent to saying that, a.s.-$P$, the measure $\rho_{N}$ converges weakly to the law of $(Z_{k_1},\ldots,Z_{k_m})$, as $N\rightarrow\infty $. For instance, by specializing (\ref{Main-ASCLT}) to the case $m=1$ one obtains that, a.s.-$P$,
\[
\frac{1}{\log N}\sum_{n=1}^N
\frac{1}{n}{\bf 1}_{\{{\rm Tr}(A^{k}_{n})-E[{\rm Tr}(A^{k}_{n})] \leq x \}}\,
{\longrightarrow} \, P\big[Z_k\leq x \big],
\]
as $N\rightarrow\infty $, for every integer $k\geq 1$ and every real $x$.
\item The content of Theorem \ref{thm-ASCLT} should be compared with the following well-known ASCLT for usual partial sums. {\sl Let $(X_{n})_{n\geq 1}$ be a sequence of real-valued independent identically distributed random
variables with
$E[X_{n}]=0$ and $E[X_{n}^{2}]=1$, and write }
$
S_{n}= \frac1{\sqrt{n}}\sum_{k=1}^n  X_k.
$
{\sl Then, almost surely, for any bounded and continuous function $\varphi:\R\to\R$,}
\begin{equation}
\label{ASCLTh}
\frac{1}{\log N}\sum_{n=1}^{N}\frac{1}{n}\varphi(S_{n}) \longrightarrow
E[\varphi(G)],\quad\mbox{\sl as $N\to\infty$}; \quad G\sim \mathscr{N}(0,1).
\end{equation}
The asymptotic relation (\ref{ASCLTh}) was first stated by L\'evy \cite{Levy} without proof,
and then forgotten for almost fifty years.
It was then rediscovered by Brosamler \cite{Brosamler} and Schatte \cite{Schatte}
and finally proved in its present form by Lacey and Philipp \cite{LaceyPhillip}.
We refer the reader to Berkes and Cs\'aki \cite{BerkesCsaki} for a universal
ASCLT covering a large class of limit theorems for partial sums, extremes, empirical
distribution functions and local times associated with independent random variables. The paper by H\"{o}rmann \cite{hormann} contains several insights into the existing literature on the subject.
\item As demonstrated in Section \ref{S : ASCLT}, in order to prove Theorem \ref{thm-ASCLT} we shall make a substantial use of a result by Ibragimov and Lifshits \cite{IL}, providing a criterion
for ASCLTs such as (\ref{ASCLTh}), not requiring that the random variables $S_n$ have the specific form of partial sums, nor that $G$ is normally distributed. Our approach is close to the one developed by Bercu {\it et al.} \cite{BNT}, in the context of ASCLTs on the Wiener space. One should also note that \cite{BNT} only deals with ASCLTs involving sequences of {\sl single} real-valued random variables (and not vectors, as in the present paper).
\end{enumerate}
}
\end{rem}

The rest of the paper is organized as follows. In Section \ref{S : UNivers} we present the universality results proved in \cite{Noupecrei3}, in a form which is convenient for our analysis. Section \ref{section3} contains a proof of (\ref{todo}). Section \ref{section4} deals with (\ref{remainder}), whereas Section \ref{S : ASCLT} focuses on the proof of Theorem \ref{thm-ASCLT}.

\section{Main tool: universality of Wiener chaos}\label{S : UNivers}
In what follows, every random object is defined on an adequate common probability space
$(\Omega, \mathscr{F}, P)$. The symbols $E$ and `${\rm Var}$' denote, respectively, the expectation and the variance associated with $P$. Also, given a finite set $B$, we write $|B|$ to indicate the cardinality of $B$. Finally, given numerical sequences $a_N,b_N$, $N\geq 1$, we write $a_N \sim b_N$ whenever $a_N/b_N \rightarrow 1$ as $N\to\infty$.

We shall now present a series of invariance principles and central limit theorems
involving sequences of homogeneous sums. These are mainly taken from
\cite{Noupecrei3} (Theorem \ref{T : NouPeRe1}),
\cite{PTu04} (Theorem \ref{T : PecTud}) and \cite{nunugio} (Theorem \ref{T : NunuGionny}). Note that the framework of \cite{Noupecrei3} is that of random variables indexed by the set of positive integers. Since in this paper we mainly deal with random variables indexed by {\sl pairs} of integers (i.e., matrix entries) we need to restate some of the findings of \cite{Noupecrei3} in terms of random variables indexed by a general (fixed) discrete countable set $A$.

\begin{defi}[Homogeneous sums]\label{Def : HomSums}
{\rm
Fix an integer $k\geq 2$. Let ${\bf Y} = \{Y_a : a\in A\} $ be a collection of
square integrable and centered independent random variables, and let $f : A^k \rightarrow \R$ be a {\sl symmetric function vanishing on diagonals} (that is, $f(a_1,...,a_k)=0$ whenever there exists $k\neq j$ such that $a_k=a_j$), and assume that $f$ has finite support. The random variable
\begin{eqnarray}
Q_k(f,{\bf Y})&=& \sum_{a_1,...,a_k \in A}f(a_1,...,a_k)Y_{a_1}\cdot\cdot\cdot Y_{a_k}
= \sum_{\{a_1,...,a_k\}\subset A^k} \!\!\!\! k!f(a_1,...,a_k)Y_{a_1}\cdot\cdot\cdot Y_{a_k}\notag\\
\label{EQ : HOMsums}
\end{eqnarray}
is called the {\sl homogeneous sum}, of order $k$, based on $f$ and ${\bf Y}$.
Clearly, $E[Q_k(f,{\bf Y})]=0$ and also, if $E(Y_a^2)=1$ for every $a\in A$, then
\begin{equation}\label{E : var}
E[Q_k(f,{\bf Y})^2] = k!\|f\|^2_{k},
\end{equation}
where, here and for the rest of the paper, we set
\[
\|f\|^2_{k} = \sum_{a_1,...,a_k \in A}f^2(a_1,...,a_k).
\]
}
\end{defi}

\smallskip

Now let ${\bf G} =\{G_a : a\in A\}$ be a collection of i.i.d. centered Gaussian random variables with unit variance. We recall that, for every $k$ and every $f$, the random variable $Q_k(f,{\bf G})$ (defined according to (\ref{EQ : HOMsums})) is an element of the $k$th {\sl Wiener chaos} associated with ${\bf G}$. See e.g. Janson \cite{Janson} for basic definitions and results on the Gaussian Wiener chaos. The next result, proved in \cite{Noupecrei3}, shows that sequences of random variables of the type $Q_k(f,{\bf G})$ have a {\sl universal character} with respect to normal approximations. The proof of Theorem \ref{T : NouPeRe1} is based on a powerful interaction between three techniques, namely: the {\sl Stein's method} for probabilistic approximations (see e.g. \cite{chen-shao}), the {\sl Malliavin calculus of variations} (see e.g. \cite{nualartbook}), and a general Lindeberg-type invariance principle recently proved by Mossel {\it et al.} in \cite{MOO}.

\begin{theorem}[Universality of Wiener chaos, see \cite{Noupecrei3}]\label{T : NouPeRe1}
Let ${\bf G}=\{G_a : a\in A\}$ be a collection of
standard centered i.i.d. Gaussian random variables, and fix integers $m\geq 1$ and $k_1,...,k_m\geq 2$.
For every $j=1,...,m$, let $\{f^{(j)}_N : N\geq 1\}$ be a sequence of functions such that
$f_N^{(j)}: A^{k_j} \rightarrow \R$ is symmetric and vanishes on diagonals. We also suppose that, for every $j=1,...,m$, the support of $f_N^{(j)}$, denoted by ${\rm supp}(f_N^{(j)})$, is such that $|{\rm supp}(f_N^{(j)})|\rightarrow \infty$, as $N\rightarrow\infty$.
Define $Q_{k_j}(f^{(j)}_N, {\bf G})$, $N\geq 1$, according to (\ref{EQ : HOMsums}). Assume that,
for every $j=1,...m$, the following sequence of variances is bounded:
\begin{equation}\label{Eq : covintro}
E[Q_{k_j}(f^{(j)}_N, {\bf G})^2] , \quad N\geq 1.
\end{equation}
Let $V$ be a $m\times m$ non-negative symmetric matrix, and let $\mathscr{N}_m(0,V)$
indicate a $m$-dimensional centered Gaussian vector with covariance matrix $V$. Then, as $N\rightarrow \infty$, the following two conditions are equivalent.
\begin{itemize}
\item[\rm (1)] The vector $\{Q_{k_j}(f^{(j)}_N, {\bf G}) : j=1,...,m\}$ converges in law to
$\mathscr{N}_m(0,V)$.
\item[\rm (2)] For every sequence ${\bf X } = \{X_a : a\in A\}$ of independent centered random variables,
with unit variance and such that $\sup_a E|X_a|^{3} <\infty$, the law of the vector
$\{Q_{k_j}(f^{(j)}_N, {\bf X}) : j=1,...,m\}$ converges to the law of $\mathscr{N}_m(0,V)$.
\end{itemize}
\end{theorem}

Note that Theorem \ref{T : NouPeRe1} concerns only homogeneous sums of order $k\geq 2$: it is easily seen (see e.g. \cite[Section 1.6.1]{Noupecrei3}) that the statement is indeed false in the case $k=1$. However, if one considers sums with a specific structure (basically, verifying some Lindeberg-type condition) one can embed sums of order one into the previous statement. A particular instance of this fact is made clear in the following statement, whose proof (combining the results of \cite{Noupecrei3} with the main estimates of \cite{MOO}) is standard and therefore omitted.

\begin{prop}\label{P : invariance+1}
For $m\geq 1$, let the kernels $\{f^{(j)}_N : N\geq 1\}$, $j=1,...,m$, verify the assumptions of Theorem \ref{T : NouPeRe1}. Let $\{a_i : i\geq 1\}$ be an infinite subset of $A$, and assume that condition (1) in the statement of Theorem \ref{T : NouPeRe1} is verified. Then, for every sequence ${\bf X } = \{X_a : a\in A\}$ of independent centered random variables, with unit variance and such that $\sup_a E|X_a|^{3} <\infty$, as $N\rightarrow\infty$ the law of the vector
$\{W_N\,;\,Q_{k_j}(f^{(j)}_N, {\bf X}) : j=1,...,m\}$, where $W_N = \frac{1}{\sqrt{N}}\sum_{i=1}^N X_{a_i}$, converges to the law of $\{N_0\, ; \, N_j:j=1,...,m\}$, where $N_0\sim\mathscr{N}(0,1)$, and $(N_1,...,N_m)\sim \mathscr{N}_m(0,V)$ denotes a centered Gaussian vector with covariance $V$, and independent of $N_0$.
\end{prop}

Theorem \ref{T : NouPeRe1} and Proposition \ref{P : invariance+1} imply that, in order to prove a CLT involving vectors of homogeneous sums based on some independent sequence ${\bf X}$, it suffices to replace ${\bf X}$ with an i.i.d. Gaussian sequence ${\bf G}$. In this way, one obtains a sequence of random vectors whose components belong to a fixed Wiener chaos. We now present two results, showing that proving CLTs for this type of random variables can be a relatively easy task: indeed, one can apply some drastic simplification of the method of moments. The first statement deals with multi-dimensional CLTs and shows that, in a Gaussian Wiener chaos setting, componentwise convergence to Gaussian always implies joint convergence. See also \cite{AirMalVie} for some connections with Stokes formula.

\begin{theorem}[Multidimensional CLTs on Wiener chaos, see \cite{Noupecrei3, PTu04}]\label{T : PecTud} Let the family ${\bf G}=\{G_a : a\in A\}$ be i.i.d. centered standard Gaussian and, for $j=1,...,m$, define the sequences $Q_{k_j}(f^{(j)}_N, {\bf G})$, $N\geq 1$, as in Theorem \ref{T : NouPeRe1} (in particular, the functions $f^{(j)}_N$ verify the same assumptions as in that theorem). Suppose that, for every $i,j=1,...,m$, as $N\rightarrow\infty$
\begin{equation}\label{jaimepas-jflg-ni-gm}
E\big[Q_{k_i}(f^{(i)}_N, {\bf G})\times Q_{k_j}(f^{(j)}_N, {\bf G})\big]\rightarrow V(i,j),
\end{equation}
where $V$ is a $m\times m$ covariance matrix.
Finally, assume that $W_N$, $N\geq 1$, is a sequence of $\mathscr{N}(0,1)$ random variables with the representation
$$
W_N = \sum_{a\in A} w_N(a)\times G_a,
$$
where the weights $w_N(a)$ are zero for all but a finite number of indices $a$, and $\sum_{a\in A} w_N(a)^2 =1$.
Then, the following are equivalent:
\begin{itemize}
\item[\rm (1)] The random vector $\{W_N\, ; \, Q_{k_j}(f^{(j)}_N, {\bf G}) : j=1,...,m\}$ converges in law to
$\{N_0\, ; \, N_j:j=1,...,m\}$, where $N_0\sim\mathscr{N}(0,1)$, and $(N_1,...,N_m)\sim \mathscr{N}_m(0,V)$ denotes a centered Gaussian vector with covariance $V$, and independent of $N_0$.
\item[\rm (2)] For every fixed $j=1,...,m$, the sequence
$Q_{k_j}(f^{(j)}_N, {\bf G})$, $N\geq 1$, converges in law
to $Z\sim\mathscr{N}\big(0,V(j,j)\big)$, that is,  to a centered Gaussian random variable with variance $V(j,j)$.
\end{itemize}
\end{theorem}

The previous statement implies that, in order to prove CLTs for vectors of homogeneous sums, one can focus on the componentwise convergence of their (Gaussian) Wiener chaos counterpart. The forthcoming Theorem \ref{T : NunuGionny} shows that this type of one-dimensional convergence can be studied by focussing exclusively on fourth moments. To put this result into full use, we need some further definitions.

\begin{defi}\label{def-contr}
{\rm Fix $k\geq 2$. Let $f : A^k \rightarrow \R$ be a (not necessarily symmetric) function vanishing on diagonals and with finite support. For every $r= 0,...,k$, the \textsl{contraction} $f\star_r f$ is the function on $A^{2d-2r}$ given by
\begin{eqnarray}\label{discretecontr}
&&f\! \star_r \!f(a_1,...,a_{2d-2r}) \\
&& =\!\!\! \sum_{(x_1,...,x_r)\in A^r} \!\!\!\!f(a_1,...,a_{k-r},x_1,...,x_r)f(a_{k-r+1},...,a_{2d-2r},x_1,...,x_r). \notag
\end{eqnarray}
Observe that (even when $f$ is symmetric) the contraction $f\star_r f $ is
not necessarily symmetric and not necessarily vanishes on diagonals. The canonical symmetrization of $f\star_r f$ is written $f\widetilde{\star}_r f$.
}
\end{defi}

\begin{theorem}[The simplified method of moments, see \cite{nunugio}]\label{T : NunuGionny} Fix $k\geq 2$. Let ${\bf G}=\{G_a : a\in A\}$ be an i.i.d. centered standard Gaussian family. Let $\{f_N : N\geq 1\}$ be a sequence of functions such that $f_N: A^{k} \rightarrow \R$ is symmetric and vanishes on diagonals. Suppose also that $|{\rm supp}(f_N)|\rightarrow \infty$, as $N\rightarrow\infty$. Assume that
\begin{equation}\label{biquet}
E[Q_{k}(f_N, {\bf G})^2] \rightarrow \sigma^2>0,\quad\mbox{as $N\rightarrow \infty$.}
\end{equation}
Then, the following three conditions are equivalent, as $N\rightarrow \infty$.
\begin{itemize}
\item[\rm (1)] The sequence $Q_{k}(f_N, {\bf G})$, $N\geq 1$, converges in law to $Z\sim\mathscr{N}(0,\sigma^2)$.
\item[\rm (2)] $E[Q_{k}(f_N, {\bf G})^4] \rightarrow 3\sigma^4$.
\item[\rm (3)] For every $r=1,...,k -1$, $\|f_N\star_r f_N\|_{2k-2r} \rightarrow 0$.
\end{itemize}
\end{theorem}

\medskip
Finally, we present a version of Theorem \ref{T : NouPeRe1} with bounds, that will lead to the proof of
Theorem \ref{Main}-(ii) provided in Section \ref{SS : Last Word}.

\begin{thm}[Universal bounds, see \cite{Noupecrei3}]\label{thm:bornes}
Let ${\bf X } = \{X_a : a\in A\}$ be a collection of independent centered random variables,
with unit variance and such that $\beta:=\sup_a E|X_a|^{3} <\infty$.
Fix integers $m\geq 1$, $k_m>...>k_1\geq 2$.
For every $j=1,...,m$, let $f^{(j)} : A^{k_j}\to\R$ be a symmetric function vanishing on diagonals.
Define $Q^j({\bf X}):=Q_{k_j}(f^{(j)}, {\bf X})$ according to (\ref{EQ : HOMsums}), and assume that
$E[Q^j({\bf X})^2]=1$ for all $j=1,\ldots,m$.
Also, assume that $K>0$ is given such that $\sum_{a\in A} \max_{1\leq j\leq m}
{\rm Inf}_a(f^{(j)})\leq K$, where
\[
{\rm Inf}_a(f^{(j)})=\sum_{\{a_2,\ldots,a_{k_j}\}\subset A^{k_j}} f^{(j)}(a,a_2,\ldots,a_{k_j})^2=
\frac{1}{(k_j-1)!}\sum_{a_2,\ldots,a_{k_j}\in A} f^{(j)}(a,a_2,\ldots,a_{k_j})^2.
\]
Let $\varphi:\R^m\to\R$ be a thrice differentiable function such that $\|\varphi''\|_\infty +
\|\varphi'''\|_\infty<\infty$, with $\|\varphi^{(k)}\|_\infty
=\max_{|\alpha|=k}\frac{1}{\alpha!}\sup_{z\in\R^m}|\partial^\alpha \varphi(z)|$.
Then, for $Z=(Z^1,\ldots,Z^m)\sim \mathscr{N}_m(0,I_m)$ (standard Gaussian vector on $\R^m$), we have
\begin{eqnarray*}
&&\big|E[\varphi(Q^1({\bf X}),\ldots,Q^m({\bf X}))]
-E[\varphi(Z)]\big|\leq \|\varphi''\|_\infty\left(\sum_{i=1}^m \Delta_{ii}+2\sum_{1\leq i<j\leq m}
\Delta_{ij}\right)\\
&&\hskip2cm +K\|\varphi'''\|_\infty
\left(\beta + \sqrt{\frac{8}{\pi}}\right)
\left[\sum_{j=1}^m (16\sqrt{2}\beta)^{\frac{k_j-1}{3}}k_j!\right]^3
\sqrt{\max_{1\leq j\leq m}\max_{a\in A} {\rm Inf}_a(f^{(j)})},
\end{eqnarray*}
where $\Delta_{ij}$, $1\leq i\leq j\leq m$, is given by
\begin{eqnarray*}
&&\frac{k_j}{\sqrt{2}}\sum_{r=1}^{k_j-1}(r-1)!\binom{k_i-1}{r-1}\binom{k_j-1}{r-1}
\sqrt{(k_i+k_j-2r)!}\big(
\|f^{(i)}\star_{k_i-r}f^{(i)}\|_{2r}\!+\!\|f^{(j)}\star_{k_j-r}f^{(j)}\|_{2r}
\big)\\
&&\hskip8cm + {\bf 1}_{\{k_i<k_j\}}\sqrt{k_j!\binom{k_j}{k_i}\|f^{(j)}\star_{k_j-k_i}f^{(j)}\|_{2k_i}}.
\end{eqnarray*}
\end{thm}


We finish this section by a useful result, which shows how the {\sl influence} ${\rm Inf}_a f$
of $f:A^k\to\R$ can be bounded by the norm of
the contraction of $f$ of order $k-1$:
\begin{prop}\label{lien-influence-contraction}
Let $f:A^k\to\R$ be a symmetric function vanishing on diagonals. Then
\[
(k-1)!\max_{a\in A} {\rm Inf}_a(f):=\max_{a\in A}
\sum_{a_2,\ldots,a_{k}\in A} f(a,a_2,\ldots,a_{k})^2
\leq \|f\star_{k-1} f\|_{2}.
\]
\end{prop}
{\it Proof}.
We have
\begin{eqnarray*}
\|f\star_{k-1}f\|^2_2 &=& \sum_{a,b\in A}
\left[
\sum_{a_2,\ldots,a_k\in A}f(a,a_2,\ldots,a_k)f(b,a_2,\ldots,a_k)
\right]^2\\
&\geq& \sum_{a\in A}
\left[
\sum_{a_2,\ldots,a_k\in A}f^2(a,a_2,\ldots,a_k)
\right]^2\\
&\geq& \max_{a\in A}
\left[
\sum_{a_2,\ldots,a_k\in A}f^2(a,a_2,\ldots,a_k)
\right]^2
=\left[(k-1)!\max_{a\in A}{\rm Inf}_a(f)\right]^2.
\end{eqnarray*}
\fin

As a consequence of Theorem \ref{thm:bornes} and Proposition \ref{lien-influence-contraction},
we immediately get the following result.
\begin{corollary}\label{cocominet}
Let ${\bf X } = \{X_a : a\in A\}$ be a collection of independent centered random variables,
with unit variance and such that $\beta:=\sup_a E|X_a|^{3} <\infty$.
Fix integers $m\geq 1$, $k_m>...>k_1\geq 1$.
For every $j=1,...,m$, let $\{f^{(j)}_N : N\geq 1\}$ be a sequence of functions such that
$f_N^{(j)}: A^{k_j} \rightarrow \R$ is symmetric and vanishes on diagonals.
Define $Q_N^j({\bf X}):=Q_{k_j}(f_N^{(j)}, {\bf X})$ according to (\ref{EQ : HOMsums}), and assume that
$E[Q_N^j({\bf X})^2]=1$ for all $j=1,\ldots,m$ and $N\geq 1$.
Let $\varphi:\R^m\to\R$ be a thrice differentiable function such that $\|\varphi''\|_\infty +
\|\varphi'''\|_\infty<\infty$.
If, for some $\alpha>0$,  $\|f_N^{(j)}\star_{k_j-r} f_N^{(j)}\|_{2r} =O(
N^{-\alpha}
)$ for all $j=1,\ldots,m$ and $r=1,\ldots, k_j-1$,
then, by noting $(Z^1,\ldots,Z^m)$ a centered Gaussian vector
such that $E[Z^iZ^j]=0$ if $i\neq j$ and $E[(Z^j)^2]=1$, we have
\[
\big|E[\varphi(Q_N^1({\bf X}),\ldots,Q_N^m({\bf X}))]
-E[\varphi(Z^1,\ldots,Z^m)]\big| =O(N^{-\alpha/2}).
\]
\end{corollary}

\section{Gaussian fluctuations of non-diagonal trace components}\label{section3}
Our aim in this section is to prove the multidimensional CLT (\ref{todo}), by using the universality results presented in Section \ref{S : UNivers}. To do this, we shall use an auxiliary collection ${\bf G}=\{G_{ij}:\,i,j\geq 1\}$
of i.i.d. copies of a $\mathscr{N}(0,1)$ random variable.

As in Section \ref{strategy}, for a given integer $k\geq 2$, we write $D_N^{(k)}$ to indicate the set of
vectors ${\bf i}=(i_1,\ldots,i_k)\in[N]^k$ such that
all the elements $(i_a,i_{a+1})$, $a=1,\ldots,k$, are different in pairs
(with the convention that $i_{k+1}=i_1$). We have the following preliminary result:

\begin{prop}\label{1dim}
For any fixed integer $k\geq 2$,
$$
N^{-k/2}\sum_{{\bf i}\in D_N^{(k)}} G_{i_1i_2}\ldots G_{i_ki_1}
\overset{{\rm Law}}{\longrightarrow}Z_k\sim \mathscr{N}(0,k)\quad\mbox{as $N\to\infty$.}
$$
\end{prop}
\begin{remark}\label{rk-verydeep}
{\rm
When $k=1$, the conclusion of the above proposition continues to be true, since
in this case we obviously have
$$
N^{-1/2}\sum_{i=1}^N G_{ii} \sim\mathscr{N}(0,1).
$$
}
\end{remark}

\medskip

\noindent{\it Proof of Proposition \ref{1dim}}: The main idea is to use the results of Section \ref{S : UNivers}, in the special case $A = \mathbb{N}^2$, that is, $A$ is the collection of all pairs $(i,j)$ such that $i,j\geq 1$. Observe that
$$
N^{-k/2}\sum_{{\bf i}\in D_N^{(k)}} G_{i_1i_2}\ldots G_{i_ki_1}
=Q_k(f_{k,N},{\bf G}),
$$
with $f_{k,N}:([N]^2)^k\to\R$ the symmetric function defined by
\begin{equation}\label{zidane}
f_{k,N}=
\frac{1}{k!}
\sum_{\sigma\in\mathfrak{S}_k}
f_{k,N}^{(\sigma)},
\end{equation}
where we used the notation
\begin{equation}\label{rrumba}
f_{k,N}^{(\sigma)}\big((a_1,b_1),\ldots,(a_k,b_k)\big)=N^{-k/2}\sum_{{\bf i}\in D_N^{(k)}}
{\bf 1}_{\{i_{\sigma(1)}=a_1,\,i_{\sigma(1)+1}=b_1\}}
\ldots
{\bf 1}_{\{i_{\sigma(k)}=a_k,\,i_{\sigma(k)+1}=b_k\}},
\end{equation}
and $\mathfrak{S}_k$ denotes the set of all permutations of $[k]$. Hence, by virtue of Theorem \ref{T : NunuGionny}, to prove Proposition \ref{1dim} it is sufficient to accomplish the following two steps:
({\it Step 1}) prove that property (3) (with $f_{k,N}$ replacing $f_N$) in the statement of Theorem \ref{T : NunuGionny} takes place, and ({\it Step 2}) show that relation (\ref{biquet}) (with $f_{k,N}$ replacing $f_N$) is verified.\\
\\
{\it Step 1}.
Let $r\in\{1,\ldots,k-1\}$. For $\sigma,\tau\in\mathfrak{S}_k$, we compute
\begin{eqnarray}\label{belleEgalite}
&&f_{k,N}^{(\sigma)}\star_r f_{k,N}^{(\tau)}\big((x_1,y_1),\ldots,(x_{2k-2r},y_{2k-2r})\big)\\
\notag
&=& N^{-k}\sum_{{\bf i},{\bf j}\in D_N^{(k)}}
{\bf 1}_{\{
i_{\sigma(1)}=x_1,\,i_{\sigma(1)+1}=y_1
\}}
\ldots
{\bf 1}_{\{i_{\sigma(k-r)}=x_{k-r},\,i_{\sigma(k-r)+1}=y_{k-r}\}}\\
&&\hskip1.5cm\times
{\bf 1}_{\{j_{\tau(1)}=x_{k-r+1},\,j_{\tau(1)+1}=y_{k-r+1}\}} \notag
\ldots
{\bf 1}_{\{j_{\tau(k-r)}=x_{2k-2r},\,j_{\tau(k-r)+1}=y_{2k-2r}\}}\\
&&\hskip1.5cm\times
{\bf 1}_{\{i_{\sigma(k-r+1)}=j_{\tau(k-r+1)},
\,i_{\sigma(k-r+1)+1}=j_{\tau(k-r+1)+1}\}}
\ldots
{\bf 1}_{\{i_{\sigma(k)}=j_{\tau(k)},
\,i_{\sigma(k)+1}=j_{\tau(k)+1}}\}. \notag
\end{eqnarray}
We now want to assess the quantity $\|f_{k,N}^{(\sigma)}\star_r f_{k,N}^{(\tau)}\|^2_{2k-2r}$. To do this, we exploit the representation (\ref{belleEgalite}) in order to write such a squared norm as a sum over $([N]^{k})^4$: as a consequence, one deduces that
$
\|f_{k,N}^{(\sigma)}\star_r f_{k,N}^{(\tau)}\|^2_{2k-2r} \leq
|F_N^{(r,\sigma,\tau)}|\,N^{-2k}
$
where $F_N^{(r,\sigma,\tau)}$ is the subset of $([N]^{k})^4$
composed of those quadruplets $({\bf i},{\bf j},{\bf a},{\bf b})$ such that
\begin{eqnarray}
&& i_{\sigma(1)}=a_{\sigma(1)},\quad
i_{\sigma(1)+1}=a_{\sigma(1)+1},\,
\ldots,\quad
i_{\sigma(k-r)}=a_{\sigma(k-r)},\quad
i_{\sigma(k-r)+1}=a_{\sigma(k-r)+1}\notag\\
&& j_{\tau(1)}=b_{\tau(1)},\quad
j_{\tau(1)+1}=b_{\tau(1)+1},\,
\ldots,\quad
j_{\tau(k-r)}=b_{\tau(k-r)},\quad
j_{\tau(k-r)+1}=b_{\tau(k-r)+1}\notag\\
&& i_{\sigma(k-r+1)}=j_{\tau(k-r+1)},\quad
i_{\sigma(k-r+1)+1}=j_{\tau(k-r+1)+1},\,
\ldots,\quad
i_{\sigma(k)}=j_{\tau(k)},\quad
i_{\sigma(k)+1}=j_{\tau(k)+1}\notag\\
&& a_{\sigma(k-r+1)}=b_{\tau(k-r+1)},\quad
a_{\sigma(k-r+1)+1}=b_{\tau(k-r+1)+1},
\ldots,\quad
a_{\sigma(k)}=b_{\tau(k)},\quad
a_{\sigma(k)+1}=b_{\tau(k)+1}.\notag\\
\label{equalities}
\end{eqnarray}
It is immediate that, among the equalities in (\ref{equalities}), the $2k$ equalities appearing in the forthcoming display (\ref{equalities2})
are pairwise disjoint (that is, an index appearing in one of the equalities does not enter into the others):
\begin{eqnarray}
&& i_{\sigma(1)}=a_{\sigma(1)},\,
\ldots,\,
i_{\sigma(k-r)}=a_{\sigma(k-r)},\quad
j_{\tau(1)}=b_{\tau(1)},\,
\ldots,\,
j_{\tau(k-r)}=b_{\tau(k-r)}
\notag\\
&& i_{\sigma(k-r+1)}= j_{\tau(k-r+1)},\,
\ldots,\quad
i_{\sigma(k)}=j_{\tau(k)},\quad
a_{\sigma(k-r+1)}=b_{\tau(k-r+1)},\,
\ldots,\,
a_{\sigma(k)}=b_{\tau(k)}.\notag\\
\label{equalities2}
\end{eqnarray}
Hence, the cardinality of $F_N^{(r,\sigma,\tau)}$ is less
than $N^{2k}$, from which we infer that $\|f_{k,N}^{(\sigma)}\star_r
f_{k,N}^{(\tau)}\|_{2k-2r}^2$ is bounded by 1. This is not sufficient for our purposes,
since we need to show that
$\|f_{k,N}^{(\sigma)}\star_r
f_{k,N}^{(\tau)}\|_{2k-2r}^2$ tends to zero as $N\to\infty$.
To prove this, it is sufficient to extract from (\ref{equalities}) one supplementary equality
which is not already written
in (\ref{equalities2}). We shall prove that this equality exists by contradiction.
Set $L=\{\sigma(s):\,1\leq s\leq k-r\}$ and $R=\{\sigma(s)+1:\,1\leq s\leq k-r\}$
(with the convention that $k+1=1$). Now assume that $R=L$. Then
$\sigma(1)+1\in R$ also belongs to $L$, so that $\sigma(1)+2\in R$. By repeating this argument,
we get that $L=R=[k]$, which is a contradiction because $r\geq 1$.
Hence, $R\neq L$. In particular, the display (\ref{equalities})
implies at least one relation involving two indices that are not already coupled in
(\ref{equalities2}).
This yields that the cardinality of $F_N^{(r,\sigma,\tau)}$ is at most $ N^{2k-1}$, and consequently
that $\|f_{k,N}^{(\sigma)}\star_r
f_{k,N}^{(\tau)}\|^2_{2k-2r}\leq N^{-1}$.
This fact implies immediately that the norms
$\|f_{k,N}\star_rf_{k,N}\|_{2k-2r}$, $r=1,\ldots,k-1$, verify
\begin{equation}\label{speed}
\|f_{k,N}\star_rf_{k,N}\|_{2k-2r} = O(N^{-1/2}),
\end{equation}
and tend to zero as $N\to\infty$.
In other words, we have proved that condition (3) in the statement of Theorem \ref{T : NunuGionny} is met.

\smallskip

\noindent {\it Step 2}. We have
$$
{\rm Var}\left(
N^{-k/2}\sum_{{\bf i}\in D_N^{(k)}}G_{i_1i_2}\ldots G_{i_ki_1}
\right)
=N^{-k}\sum_{{\bf i},{\bf j}\in D_N^{(k)}}
E[G_{i_1i_2}\ldots G_{i_ki_1}
G_{j_1j_2}\ldots G_{j_kj_1}].
$$
For fixed ${\bf i},{\bf j}\in D_N^{(k)}$,
observe that the expectation $E[G_{i_1i_2}\ldots G_{i_ki_1}
G_{j_1j_2}\ldots G_{j_kj_1}]$ can only be zero or one. Moreover, it is one
if and only if, for all $s\in [k]$,
there is exactly one $t\in[k]$ such that $(i_s,i_{s+1})=(j_t,j_{t+1})$.
In this case, we define $\sigma\in\mathfrak{S}_k$ as
the bijection of $[k]$ into itself which maps each $s$ to the corresponding $t$ and we have, for all $s\in[k]$,
\begin{equation}\label{is}
i_s=j_{\sigma(s)}=j_{\sigma(s-1)+1}.
\end{equation}
To summarize, one has
that
${\rm Var}\left(
N^{-k/2}\sum_{{\bf i}\in D_N^{(k)}}G_{i_1i_2}\ldots G_{i_ki_1}
\right)$
equals
\begin{equation}\label{variance}
N^{-k}\sum_{\sigma\in\mathfrak{S}_k}\big|
\big\{ ({\bf i},{\bf j})\in (D_N^{(k)})^2:\,(i_s,i_{s+1})=(j_{\sigma(s)},j_{\sigma(s)+1})\,\,\,
\mbox{for all $s\in[k]$}\big\}
\big|.
\end{equation}
If $\sigma\in\mathfrak{S}_k$ is such that $\sigma(s)=\sigma(s-1)+1$ for all $s$
(it is easily seen that there are exactly $k$ permutations verifying this property in $\mathfrak{S}_k$),
we get $k$ different conditions by letting $s$ run over $[k]$
in (\ref{is}), so that
$$
\big\{ ({\bf i},{\bf j})\in (D_N^{(k)})^2:\,(i_s,i_{s+1})=(j_{\sigma(s)},j_{\sigma(s)+1})\,\,\,
\mbox{for all $s\in[k]$}\big\}= N^k+O(N^{k-1}),\quad \mbox{as $N\to\infty$}.
$$
In contrast, if $\sigma\in\mathfrak{S}_k$ is {\sl not} such that $\sigma(s)=\sigma(s-1)+1$ for all $s$,
then by letting $s$ run over $[k]$, one deduces from (\ref{is}) at least $k+1$ different conditions,
so that, in this case,
$$
\big\{ ({\bf i},{\bf j})\in (D_N^{(k)})^2:\,(i_s,i_{s+1})=(j_{\sigma(s)},j_{\sigma(s)+1})\,\,\,
\mbox{for all $s\in[k]$}\big\}=O(N^{k-1}),\quad \mbox{as $N\to\infty$}.
$$
Taking into account these two properties together with the representation (\ref{variance}), we deduce that
the variance of \[
N^{-k/2}\sum_{{\bf i}\in D_N^{(k)}}G_{i_1i_2}\ldots G_{i_ki_1}
\] tends to $k$
 as $N\to\infty$. It follows that the required property (\ref{biquet}) in Theorem \ref{T : NunuGionny} (with $\sigma^2=k$) is met.

\smallskip

The proof of Proposition \ref{1dim} is concluded.
\fin

\begin{rem}{\rm
By inspection of the previous proof, one also deduces that, for every $k\geq 2$,
there exists a constant $C_k$ (independent of $N$) such that, for all $N\geq 1$,
\begin{equation}\label{prop1}
\left|{\rm Var}\left(N^{-k/2}\sum_{{\bf i}\in D_N^{(k)}}G_{i_1i_2}\ldots G_{i_ki_1}\right)- k\right|\leq\frac{C_k}{N}.
\end{equation}
}
\end{rem}

The {\sl multidimensional} version of Proposition \ref{1dim} reads as follows:
\begin{prop}\label{P : multiNonD}
Fix $m\geq 1$, as well as integers $k_m> \ldots> k_1\geq 2$.
Then, as $N\to\infty$,
\begin{eqnarray}\label{buxethude}
&&\left(N^{-1/2}\sum_{i=1}^N G_{ii},\,\,\,
N^{-\frac{k_1}{2}}  \sum_{{\bf i}\in D_N^{(k_1)}} G_{i_1 i_2}\cdot\cdot\cdot G_{i_{k_1} i_1}
,\,\ldots \right. \\
&&\quad\quad\quad\quad\quad\quad\quad\quad \quad \left. \ldots,\,
N^{-\frac{k_m}{2}}  \sum_{{\bf i}\in D_N^{(k_m)}} G_{i_1 i_2}\cdot\cdot\cdot G_{i_{k_m} i_1}
\right)
\stackrel{\rm Law}{\longrightarrow} \big(Z_1,Z_{k_1} ,..., Z_{k_m}\big),\notag
\end{eqnarray}
where ${\bf Z}=\{Z_k:k\geq 1\}$ denotes a collection of independent centered Gaussian random
variables such that, for every $k\geq 1$, $E(Z_k^2)=k$.
\end{prop}
{\it Proof}:
It is an application of Theorem \ref{T : PecTud}, in the following special case:
\begin{itemize}
\item[ - ] $w_N(i,j) = \frac{1}{\sqrt{N}}$, if $i=j\leq N$ and $w_N(i,j) = 0$ otherwise;
\item[ - ] $V$ is equal to the diagonal matrix such that $V(a,b)=0$ if $a\neq b$ and $V(a,a) =k_a$, for $a=1,...,m$;
\item[ - ] for $j=1,...,m$, $f_N^{(j)} = f_{k_j,N}$, where we used the notation (\ref{zidane}).
\end{itemize}
Indeed, in view of Proposition \ref{1dim}, one has that condition (2) in the statement of Theorem \ref{T : PecTud} is satisfied.
Moreover, for fixed $a\neq b$ and since ${\bf G}$ consists of a collection of independent and
centered (Gaussian) random variables, it is clear
that, for all $N$,
$$
E\left[
\sum_{{\bf i}\in D_N^{(k_a)}} G_{i_1i_2}\ldots G_{i_{k_a}i_1}
\times
\sum_{{\bf j}\in D_N^{(k_b)}} G_{j_1j_2}\ldots G_{j_{k_b}j_1}
\right]=0,
$$
so that condition (\ref{jaimepas-jflg-ni-gm}) is met. The proof is concluded.
\fin

\medskip

By combining Proposition \ref{P : multiNonD} and Proposition \ref{P : invariance+1}, we can finally deduce the following general result for non-diagonal trace components.

\begin{cor}\label{P : megaP} For $N\geq 2$, let $X_N$ be the $N\times N$ random matrix given by (\ref{Def : RM}), where the reference random variable $X$ has mean zero, unit variance and finite absolute third moment. Fix $m\geq 1$, as well as integers $2\leq k_1<\ldots< k_m$. Then, the CLT (\ref{todo}) takes place,
with ${\bf Z}=\{Z_k:k\geq 1\}$ denoting a sequence of independent centered Gaussian random
variables such that, for every $k\geq 1$, $E(Z_k^2)=k$.
\end{cor}

\begin{rem}{\rm In order to prove Corollary \ref{P : megaP}, one only needs the existence of third moments. Note that, as will become clear in the following Section \ref{section4}, moments of higher orders are necessary for our proof of (\ref{remainder}).
}
\end{rem}

\section{The remainder: combinatorial bounds on partitioned chains and proof of Theorem \ref{Main}}\label{section4}

Fix an integer $k\geq 2$. From section \ref{strategy},
recall that $D_N^{(k)}$ denotes the subset
of vectors ${\bf i}=(i_1,\ldots,i_k)\in[N]^k$ such that
all the elements $(i_a,i_{a+1})$, $a=1,\ldots,k$, are different in pairs (with the convention that $i_{k+1}=i_1$).
From the Introduction, recall that $X$ is a centered random variable,
having unit variance and with finite moments of all orders.
Let also ${\bf X}=\{X_{ij}:i,j\geq 1\}$ be a collection of i.i.d. copies of $X$. In the present section, our aim is
to prove the asymptotic relation (\ref{remainder}), that is
\begin{prop}\label{remainder-prop}
For every $k\geq 2$, as $N\to\infty$,
\begin{equation}\label{Graal}
E(R_N(k)^2) = {\rm Var}\left(
N^{-k/2}\sum_{{\bf i}\not\in D_N^{(k)}}
\big[X_{i_1i_2}\ldots X_{i_ki_1}
-E(X_{i_1i_2}\ldots X_{i_ki_1})\big]
\right)=O(N^{-1}).
\end{equation}
\end{prop}
The proof of Proposition \ref{remainder-prop} is detailed in Section \ref{SS : Last Word}, and builds on several combinatorial estimates derived in Sections \ref{SS : CombDefs}--\ref{SS : CombBounds}. To ease the reading of the forthcoming material, we now provide an intuitive outline of this proof.

\smallskip

\noindent{\bf Remark on notation.} Given an integer $k\geq 2$, we denote by $\mathcal{P}(k)$ the collection of all
partitions of $[k]=\{1,...,k\}$. Recall that a partition $\pi\in\mathcal{P}(k)$ is an object of the type
$\pi = \{B_1,...,B_r\}$, where the $B_j$'s are disjoint and non-empty subsets of $[k]$, called {\sl blocks},
such that $\cup_{j=1,...,r} B_j = [k]$. Given $a,x\in[k]$ and $\pi\in\mathcal{P}(k)$, we write $a\stackrel{\pi}{\sim} x$ whenever $a$ and $x$ are in the same block of $\pi$. We also use the symbol $\hat{1}$ to indicate the one-block partition $\hat{1}=\{[k]\}$  (this is standard notation from combinatorics -- see e.g. \cite{Stanley}). In this section, for the sake of simplicity and because $k$ is fixed, we write $D_N$ instead of $D_N^{(k)}$.

\subsection{Sketch of the proof of Proposition \ref{remainder-prop}}\label{SS : Sketch}

Our starting point is the following elementary decomposition:
\[
[N]^k\setminus D_N = \bigcup_{\pi\in\mathcal{Q}(k)}A_N(\pi),
\]
where $\mathcal{Q}(k)$ stands for the collection of all partitions of $[k]$ containing at least one block of cardinality $\geq 2$,
and $A_N(\pi)$ is the collection of all vectors ${\bf i}\in[N]^k$
such that the equality $(i_a,i_{a+1})=(i_x,i_{x+1})$ holds if and only if $a\stackrel{\pi}{\sim} x$.
Using this decomposition, one sees immediately that, in order to show (\ref{Graal}), it is sufficient
to prove that, for each {\sl fixed} $\pi\in\mathcal{Q}(k)$,
the quantity
\begin{eqnarray}\label{sum}
&&   {\rm Var}\left( N^{-k/2}\sum_{{\bf i}\in A_N(\pi)}
\big[X_{i_1i_2}\ldots X_{i_ki_1}
-E(X_{i_1i_2}\ldots X_{i_ki_1})\big]  \right)    \\
&& =N^{-k}
\!\!\!\!\!\!\!\!\!\!
\sum_{({\bf i},{\bf j})\in A_N(\pi)\times A_N(\pi)}
\!\!\!\!\!\!\!\!\!\!
\left[\,E(X_{i_1i_2}\ldots X_{i_ki_1}
X_{j_1j_2}\ldots X_{j_kj_1})
-E(X_{i_1i_2}\ldots X_{i_ki_1})
E(X_{j_1j_2}\ldots X_{j_kj_1})\,\right] \notag
\end{eqnarray}
is $O(N^{-1})$, as $N\to\infty$. Let $G_N(\pi)$ denote the subset of
pairs $({\bf i},{\bf j})\in A_N(\pi)\times A_N(\pi)$ such that
the following non-vanishing condition is in order:
\begin{equation}\label{nonzero}
E(X_{i_1i_2}\ldots X_{i_ki_1}
X_{j_1j_2}\ldots X_{j_kj_1})
-E(X_{i_1i_2}\ldots X_{i_ki_1})
E(X_{j_1j_2}\ldots X_{j_kj_1})\neq 0.
\end{equation}
Hence
\begin{eqnarray}\label{sumbis}
&&   {\rm Var}\left( N^{-k/2}\sum_{{\bf i}\in A_N(\pi)}
\big[X_{i_1i_2}\ldots X_{i_ki_1}
-E(X_{i_1i_2}\ldots X_{i_ki_1})\big]  \right)    \\
&& =N^{-k}
\!\!\!\!\!
\sum_{({\bf i},{\bf j})\in G_N(\pi)}
\!\!\!\!\!
\left[\,E(X_{i_1i_2}\ldots X_{i_ki_1}
X_{j_1j_2}\ldots X_{j_kj_1})
-E(X_{i_1i_2}\ldots X_{i_ki_1})
E(X_{j_1j_2}\ldots X_{j_kj_1})\,\right]. \notag
\end{eqnarray}
Due to the finite moment assumptions for $X$, and
by appling the generalized H\"older inequality,
it is clear that, for a generic pair $({\bf i},{\bf j})$,
\[
\big|E(X_{i_1i_2}\ldots X_{i_ki_1}
X_{j_1j_2}\ldots X_{j_kj_1})
-E(X_{i_1i_2}\ldots X_{i_ki_1})
E(X_{j_1j_2}\ldots X_{j_kj_1}) \big|
\leq 2\,E(|X|^{2k})<\infty.
\]
It follows that, in order to prove that the sum in (\ref{sumbis}) is $O(N^{-1})$, it is enough to show that
\begin{equation}\label{gn}
\big|G_N(\pi)\big|\leq\Theta(k,\pi) N^{k-1},
\end{equation}
for some constant $\Theta(k,\pi)$ not depending on $N$.
Our way of proving (\ref{gn}) is to show that, if $({\bf i},{\bf j})$ denotes a {\sl generic} element
of $G_N(\pi)$, then, necessarily, there
exists at least $k+1$ equalities between the $2k$ indices $i_1,\ldots,i_k,$ $j_1,\ldots,j_k$
of $({\bf i},{\bf j})$. Note that by `equality'
we just mean the existence of two {\sl different} integers $a,b\in[k]$ such that
$i_a=i_b$ or $j_a=j_b$, or the existence of two integers $a,b\in[k]$ such that
$i_a=j_b$. Proving this fact implies that the $2k$ indices of a generic elements
$({\bf i},{\bf j})$ of $G_N(\pi)$
have at most $k-1$ {\sl degrees of freedom} (see Point {\bf 7} of Section 4.2
for a precise definition), so that (\ref{gn}) holds immediately --- the constant $\Theta(k,\pi)$ merely counting the number of ways in which the $k+1$ equalities can be consistently distributed among the indices composing $({\bf i},{\bf j})$. In order to extract these $k+1$ equalities between the $2k$ indices of a generic
element $({\bf i},{\bf j})$ of $G_N(\pi)$, we will consider two cases, according as the partition
$\pi\in\mathcal{Q}(k)$
contains at least one singleton or not.\\

\noindent {\it Case A: No singletons in $\pi$}. By definition of $A_N(\pi)$,
and due to the absence of singleton in $\pi$, we already see that there
are at least $k/2$ or $(k+1)/2$ (according to the evenness of $k$) equalities between
the $k$ indices of ${\bf i}$ (resp. ${\bf j}$).
Moreover,
the non-vanishing condition (\ref{nonzero}) implies that there is at least one further equality between one index of ${\bf i}$
and one index of ${\bf j}$. So, we proved the existence of $k+1$ equalities between
the $2k$ indices of $({\bf i},{\bf j})$, and the proof of (\ref{gn}) in the Case A is done.\\
\\
{\it Case B: At least one singleton in $\pi$}. Let $S$ denote the
collection of the singleton(s) of $\pi$.
In order for (\ref{nonzero}) to be true, observe that, for all $s\in S$,
we must have $(j_s,j_{s+1})=(i_a,i_{a+1})$ for some $a\in[k]$. In particular,
this means that there exist $|S|$ equalities of the type $j_s=i_a$ for the indices composing $({\bf i},{\bf j})$.
Also, by definition of the objects we are dealing with, for all $t\in [k]\setminus S$,
we must have $(i_t,i_{t+1})
=(i_a,i_{a+1})$ for some $a$, different from $t$, in the same $\pi$-block as $t$. Of course, the same must hold with
$i$ replaced by $j$.
Hence,
in order for (\ref{gn}) to be true, it remains to produce one equality between indices
that has not been already considered.
We mentioned above that
for all $t\in[k]\setminus S$, there exists
$a$, different from $t$ and in the same block as $t$,
such that
$j_t=j_a$. Hence, to conclude it remains to
show that we have $j_t=j_a$ for at least one integer $t$ belonging to
$[k]\setminus S$
and one integer $a$ {\sl not} belonging to the same block as $t$. Since, by assumption,
$\pi$ contains
at least one singleton
and one block of cardinality $\geq 2$ (indeed, $\pi\in \mathcal{Q}(k)$), without loss
of generality (up to relabeling the indices according to a cyclic permutation of $[k]$), we can assume that $S$ contains the singleton $\{k\}$.
Consider now the singleton $\{s^*\}$ of $S$, where $s^*$ is defined as the greatest of the integers $m$ such that
$\{m\}$ is adjacent from the right to a block, say $B_{u*}$, of cardinality $\geq 2$. For a particular example of this situation,
see the diagram in Fig. \ref{fig1}, where each row represents the same partition of $[7]$ having $s^*=6$ (see Point {\bf 3}. in the subsequent Section \ref{SS : CombDefs} for a formal construction of diagrams). To finish the proof, once again we split it into two cases:\\
\\
{\it Case B1: The block $B_{u*}$ contains two consecutive integers}. This assumption
implies that $j_x=j_t=j_{t+1}$ for all $x,t\in B_{u*}$. Since $\{a\}$ is adjacent from the right
to $B_{u*}$, we have $j_a=j_t$ for all $t\in B_{u*}$, which is exactly what we wanted to show.\\
\\
{\it Case B2: The block $B_{u*}$ does not contain two consecutive integers}.
Fig. \ref{fig8} is an illustrative example of such situation, where each row represents the same partition of $[8]$, with $s^*=7$. As we see
on this picture, we have necessarily $j_7 = j_5$, yielding the desired additional equality, which could not be extracted from the previous discussion. In Section \ref{SS : CombBounds}, it is shown that this line of reasoning can be extended to general situations.

\begin{rem}{\rm
The sketch given above contains all the main ideas entering in the proof of Proposition \ref{remainder-prop}.
The reader not interested in technical combinatorial details, can then go directly to Section
\ref{SS : Last Word}, where the proof of Theorem \ref{Main} is concluded. The subsequent Sections
\ref{SS : CombDefs}--\ref{SS : CombBounds} fill the gaps of the above sketch,
by providing exact definitions as well as complete formal arguments leading to the estimate (\ref{Graal}).
}
\end{rem}

\subsection{Definitions}\label{SS : CombDefs}

In the following list, we introduce some further definitions that are needed for the analysis developed in the rest of this
section.\\
\\
{\bf 1.} Fix integers $N,k \geq 2$. A {\sl chain} $c$ of length $2k$, built from $[N]$, is an object given by the juxtaposition of $2k$ pairs of integers of the type
\begin{equation}\label{chain}
c = (i_1,i_2)(i_2,i_3)...(i_k,i_1)(j_1,j_2)(j_2,j_3)...(j_k,j_1),
\end{equation}
where $i_a,j_x \in [N]$, for $a,x=1,...,k$. The class of all chains of length $2k$ built from $[N]$ is denoted by $C(2k,N)$. As a notational convention, we will use the letter $i$ to write the first $k$ pairs in the chain, and the letter $j$ to write the remaining ones. For instance, an element of $C(6,5)$ (that is, a chain of length 6 built from the set $\{1,2,3,4,5\}$) is
\[
c = (1,5)(5,1)(1,1)(3,3)(3,3)(3,3),
\]
where $i_1=1$, $i_2=5$, $i_3=1$, $j_1=j_2=j_3=3$. According to the graphical conventions given below (at Point 3 of the present list) we will sometimes say that $(i_1,i_2)(i_2,i_3)...(i_k,i_1)$ and $(j_1,j_2)(j_2,j_3)$ $...(j_k,j_1)$ are, respectively, the {\sl upper sub-chain} and the {\sl lower sub-chain} associated with the chain $c$ in (\ref{chain}). For instance, in the previous example the upper sub-chain is $(1,5)(5,1)(1,1)$, whereas the lower one is $(3,3)(3,3)(3,3)$. We shall say that $(i_l,i_{l+1})$ is the $l$th pair in the upper sub-chain of $c$ (and similarly for the elements of the lower sub-chain). We shall sometimes call $i_a$ the {\sl left index} of the pair $(i_a,i_{a+1})$. Also, we use the convention $i_{k+1}=i_1$ and $j_{k+1}=j_1$. Of course, a chain is completely determined by the left indices of its pairs.\\
\\
{\bf 2.} Let $\pi  \in \mathcal{P}(k)$ be a partition of $[k]$. Recall that, for $ a,b\in [k]$, we write $a \stackrel{\pi}{\sim} b$ to indicate that $a$ and $b$ belong to the same block of $\pi$. We say that a chain $c$ as in (\ref{chain}) {\sl has partition} $\pi$ if, for every $a,b\in [k]$, the following double implications take place: (i) $(i_a,i_{a+1})=(i_b,i_{b+1})$ if and only if $a \stackrel{\pi}{\sim} b$, and (ii) $(j_a,j_{a+1})=(j_b,j_{b+1})$ if and only if $a \stackrel{\pi}{\sim} b$. In other words, a chain has partition $\pi$ if and only if the partitions of $[k]$ induced by the identical pairs in its upper and lower sub-chain are both equal to $\pi$, that is (with the notation of Section \ref{SS : Sketch}), if and only if $(i_1,...,i_k),(j_1,...,j_k)\in A_N(\pi)$. For instance, take $k=4$ and $\pi = \{\{1,3\},\{2,4\}\}$. Then, the following chain built from $[3]$ has partition $\pi$:
    \[
    c=(1,2)(2,1)(1,2)(2,1)(3,1)(1,3)(3,1)(1,3).
    \]
    Note the `only if' part in the definition given above, implying that, if a chain has partition $\pi$ and if $x$ and $y$ are not in the same block of $\pi$, then necessarily $(i_x,i_{x+1}) \neq (i_y,i_{y+1}) $ and $(j_x,j_{x+1}) \neq (j_y,j_{y+1})$. This yields in particular that a chain cannot have two different partitions.\\
\\
{\bf 3.} Given $k\geq 2$, we shall sometimes represent a generic chain with partition $\pi\in\mathcal{P}(k)$ by means of {\sl diagrams}. These diagrams are mnemonic devices composed of an upper row and a lower row, of $k$ dots each. These rows represent, respectively, the upper and lower sub-chain of a given chain, in such a way that the $l$th dot (from left to right) in the upper (resp. lower) row corresponds the $l$th pair in the upper (resp. lower) sub-chain. Each block $B$ of the partition $\pi$ is represented by two closed curves: the first one is drawn around the dots of the upper row corresponding to the pairs $(i_a,i_{a+1})$ verifying $a\in B$; the second one is drawn around the dots of the lower row corresponding to those $(j_x,j_{x+1})$ verifying $x\in B$. The resulting diagram is the superposition of two identical combinations of dots and curves. Note that the shape of the diagram does not depend on $N$.
    For instance, the diagram in Fig. \ref{fig1} corresponds to the case $k=7$, and $\pi=\{\{1,4,5\},\{2\},\{3\},\{6\},\{7\}\}$,\footnote{A chain with partition $\pi$ as in Fig. \ref{fig1} is $$c = (1,1)(1,2)(2,1)(1,1)(1,1)(1,3)(3,1)(1,1)(1,4)(4,1)(1,1)(1,1)(1,5)(5,1).$$} whereas the diagram in Fig. \ref{fig2} corresponds to $k=6$ and the
one-block partition $\hat{1}=\{[6]\}$.
{
\begin{figure}[h]
\begin{center}
\psset{unit=1cm}

\begin{pspicture}(0,-1.04)(5.92,1.21)
\psframe[linewidth=0.02,dimen=outer](5.92,1.21)(0.0,-0.99)
\psdots[dotsize=0.12](1.4,0.67)
\psellipse[linewidth=0.02,dimen=outer](1.4,0.67)(0.2,0.2)
\psdots[dotsize=0.12](4.6,0.67)
\psellipse[linewidth=0.02,dimen=outer](4.6,0.67)(0.2,0.2)
\psdots[dotsize=0.12](0.6,0.67)
\psdots[dotsize=0.12](3.02,0.67)
\psdots[dotsize=0.12](2.2,0.67)
\psellipse[linewidth=0.02,dimen=outer](2.2,0.67)(0.2,0.2)
\psdots[dotsize=0.12](3.8,0.67)
\psbezier[linewidth=0.02](0.26,0.69)(0.26,-0.03)(3.34,-0.03)(3.34,0.69)
\psbezier[linewidth=0.02](0.78,0.71)(0.78,0.17)(2.86,0.17)(2.86,0.71)
\psarc[linewidth=0.02](0.52,0.69){0.26}{-0.0}{180.0}
\psline[linewidth=0.02cm](2.18,0.69)(2.2,0.67)
\psdots[dotsize=0.12](5.4,0.67)
\psellipse[linewidth=0.02,dimen=outer](5.4,0.67)(0.2,0.2)
\psarc[linewidth=0.02](3.43,0.44){0.63}{21.037512}{154.65382}
\psarc[linewidth=0.02](3.46,0.61){0.14}{26.565052}{151.69925}
\psarc[linewidth=0.02](3.8,0.81){0.26}{-149.03624}{-29.623749}
\psdots[dotsize=0.12](1.42,-0.33)
\psellipse[linewidth=0.02,dimen=outer](1.42,-0.33)(0.2,0.2)
\psdots[dotsize=0.12](4.62,-0.33)
\psellipse[linewidth=0.02,dimen=outer](4.62,-0.33)(0.2,0.2)
\psdots[dotsize=0.12](0.62,-0.33)
\psdots[dotsize=0.12](3.04,-0.33)
\psdots[dotsize=0.12](2.22,-0.33)
\psellipse[linewidth=0.02,dimen=outer](2.22,-0.33)(0.2,0.2)
\psdots[dotsize=0.12](3.82,-0.33)
\psbezier[linewidth=0.02](0.28,-0.31)(0.28,-1.03)(3.36,-1.03)(3.36,-0.31)
\psbezier[linewidth=0.02](0.8,-0.29)(0.8,-0.83)(2.88,-0.83)(2.88,-0.29)
\psarc[linewidth=0.02](0.54,-0.31){0.26}{-0.0}{180.0}
\psline[linewidth=0.02cm](2.2,-0.31)(2.22,-0.33)
\psdots[dotsize=0.12](5.42,-0.33)
\psellipse[linewidth=0.02,dimen=outer](5.42,-0.33)(0.2,0.2)
\psarc[linewidth=0.02](3.46,-0.57){0.64}{23.198591}{154.65382}
\psarc[linewidth=0.02](3.48,-0.39){0.14}{29.291363}{151.69925}
\psarc[linewidth=0.02](3.82,-0.19){0.26}{-149.03624}{-29.623749}
\end{pspicture}
\end{center}
\caption{a chain with a five-block partition}\label{fig1}
\end{figure}
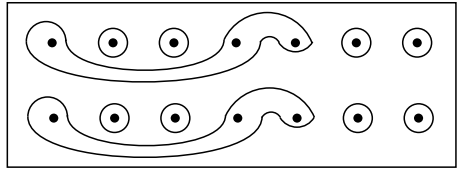
}

{
\begin{figure}[h]
\begin{center}
\psset{unit=1cm}
\begin{pspicture}(0,-0.84)(5.06,0.84)
\psdots[dotsize=0.12](0.54,0.4)
\psdots[dotsize=0.12](1.32,0.4)
\psdots[dotsize=0.12](2.12,0.4)
\psdots[dotsize=0.12](2.92,0.4)
\psdots[dotsize=0.12](3.7,0.4)
\psdots[dotsize=0.12](4.5,0.4)
\psframe[linewidth=0.02,dimen=outer](5.06,0.84)(0.0,-0.84)
\psellipse[linewidth=0.02,dimen=outer](2.54,0.39)(2.32,0.33)
\psdots[dotsize=0.12](0.54,-0.38)
\psdots[dotsize=0.12](1.32,-0.38)
\psdots[dotsize=0.12](2.12,-0.38)
\psdots[dotsize=0.12](2.92,-0.38)
\psdots[dotsize=0.12](3.7,-0.38)
\psdots[dotsize=0.12](4.5,-0.38)
\psellipse[linewidth=0.02,dimen=outer](2.54,-0.39)(2.32,0.33)
\end{pspicture}
\end{center}
\caption{a chain with a one-block partition}\label{fig2}
\end{figure}
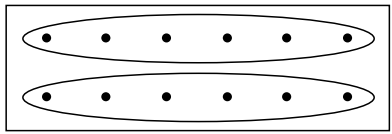
}
\noindent
{\bf 4.} In general, given a chain $c$ as in (\ref{chain}) with partition $\pi=\{B_1,...,B_r\}$ as at Point 2 of the present list, we shall say the the block $B_u$ of the upper sub-chain {\sl corresponds} to the block $B_v$ of the lower sub-chain, whenever $(i_a,i_{a+1}) = (j_x,j_{x+1})$ for every $a\in B_u $ and every $x \in B_v$. Note that one given block $B_u$ in the upper sub-chain cannot correspond to more than one block in the lower sub-chain. For $\pi = \{B_1,...,B_r\} \in \mathcal{P}(k)$, we shall now define a class of chains $C_\pi (2k,N) \subset C(2k,N)$, whose elements have partition $\pi$ and are characterized by two facts: the associated upper and lower sub-chains have at least one pair in common, and ``no singletons are left on their own''. Formally, the class $C_\pi (2k,N)$ is defined as follows (recall that we use the letter $i$ for the elements of the upper sub-chain, and the letter $j$ for the elements of the lower sub-chain). (i) If $|B_t|\geq 2$ for every $t=1,...,r$, then $C_\pi (2k,N)$ is the collection of all chains of partition $\pi$ verifying that there exists $a,x\in[k]$ such that the block $B_a$ in the upper sub-chain corresponds to the block $B_x$ in the lower sub-chain. (ii) If $\pi$ contains at least one singleton, then $C_\pi (2k,N)$ is the collection of all chains of partition $\pi$ such that every singleton in the upper (resp. lower) sub-chain corresponds to a block of the lower (resp. upper) subchain, that is: for every $\{ a \}\in\pi$, there exists $u=1,...,r$ such that $(i_a,i_{a+1})=(j_l,j_{l+1})$ for every $l \in B_u$, and, for every $\{ x \}\in\pi$, there exists $v=1,...,r$ such that $(j_x,j_{x+1})=(j_s,j_{s+1})$ for every $s \in B_v$. For instance, if $k=3$ and $\pi = \{[3]\}$, then one element of $C_\pi(6,5)$ is
    \[
   c = (5,5)(5,5)(5,5)(5,5)(5,5)(5,5).
    \]
If $k=6$ and $\pi = \{\{1,2,3\},\{4\},\{5\},\{6\}\}$, then one element of $C_\pi(12,5)$ is
    \[
    c = (1,1)(1,1)(1,1)(1,2)(2,5)(5,1)(2,2)(2,2)(2,2)(2,5)(5,1)(1,2).
    \]
{\bf 5.} Fix $ k,N \geq 2 $, as well as a partition $\pi =\{B_1,...,B_r\}\in \mathcal{P}(k)$. Given two subsets $U,V \subset [r]$ such that $|U|=|V|$, let $R : U \rightarrow V : u\mapsto R(u)$ be a bijection from $U$ onto $V$. We shall denote by $C_\pi ^R (2k,N)$ the subset of $C_\pi(2k,N)$ composed of those chains $c \in C_\pi(2k,N)$ such that the block $B_u$ in the upper sub-chain corresponds to the block $B_{R(u)}$ in the lower sub-chain. When $U=\{u\}$ and $V=\{v\}$ are singletons, we shall simply write $C_\pi ^{u,v} (2k,N)$ to indicate the set of those $c \in C_\pi(2k,N)$ such that the block $B_u$ in the upper sub-chain corresponds to the block $B_{v}$ in the lower sub-chain. For instance, the chain
        \[
        c_1 = (1,1)(1,1)(1,2)(2,5)(5,1)(2,2)(2,2)(2,5)(5,1)(1,2)
        \]
        is an element of $C_\pi^R(10, 4)$, where $\pi = \{B_1,B_2,B_3,B_4\}=\{\{1,2\}, \{3\}, \{4\}, \{5\}\}$, $U=V=\{2,3,4\}$, and $R(2) = 4$, $R(3)=2$ and $R(4)=3$. The chain
        \[
        c_2 = (3,3)(3,3)(3,3)(3,3)
        \]
        belongs to $C_{\hat{1}}^{1,1}(4, 3)$, where $\hat{1} =\{B_1\}= \{[2]\}$. Note that the definition of $C_\pi ^R (2k,N)$ does not give any information concerning the blocks of the upper and lower sub-chains that do not belong, respectively, to the domain and the image of $R$. In other words, for a chain $c\in C_\pi ^R (2k,N)$, one can have that the block $B_u$ in the upper sub-chain corresponds to the block $B_v$ in the lower sub-chain even if $u\in \!\!\!\!\!/ \, U$ and $v\in \!\!\!\!\!/ \, V$. For instance, the chain
        \[
        c = (1,1)(1,1)(1,2)(2,5)(5,1)(1,1)(1,1)(1,2)(2,5)(5,1)
        \]
        is counted as an element of $C_\pi^R(10, 4)$, where \[\pi = \{B_1,B_2,B_3,B_4\}=\{\{1,2\}, \{3\}, \{4\}, \{5\}\},\] $U=V=\{2,3,4\}$, and $R(u) = u$, for $u=2,3,4$.\\
\\
{\bf 6.} Fix $ k,N \geq 2 $, as well as a partition $\pi =\{B_1,...,B_r\}\in \mathcal{P}(k)$.
Given a bijection $R : U \rightarrow V $ as at Point 5 above,
we shall  represent a generic element of the class $C_\pi ^R (2k,N)$ by means of a diagram
built as follows: first (i) draw the diagram associated with the class $C_\pi (2k,N)$,
as explained at Point 3 of the present list, then (ii) for every pair of blocks $B_u$ and $B_v$
such that $u\in U$, $v\in V$ and $v=R(u)$ (note that $B_u$ is in the upper sub-chain, and $B_v$
in the lower sub-chain), draw a segment linking a representative element of $B_u$ with a
representative element of $B_v$. For instance, the class $C_\pi ^R (10,N)$, associated with the
chain $c_1$ appearing at Point 5 above, is represented by the diagram appearing in Fig. \ref{fig3}, whereas the chain $c_2$ is associated with the class $C_{\hat{1}} ^{1,1} (4,3)$, whose diagram is drawn in Fig. \ref{fig4}.
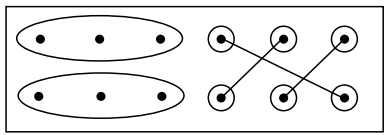
\begin{figure}[h]
\begin{center}
\psset{unit=1cm}
\begin{pspicture}(0,-0.84)(4.98,0.84)
\psdots[dotsize=0.12](0.46,0.4)
\psdots[dotsize=0.12](1.24,0.4)
\psdots[dotsize=0.12](2.04,0.4)
\psellipse[linewidth=0.02,dimen=outer](1.24,0.41)(1.1,0.31)
\psdots[dotsize=0.12](2.84,0.4)
\psdots[dotsize=0.12](0.44,-0.36)
\psdots[dotsize=0.12](1.26,-0.36)
\psdots[dotsize=0.12](2.06,-0.36)
\psellipse[linewidth=0.02,dimen=outer](1.26,-0.35)(1.1,0.31)
\psframe[linewidth=0.02,dimen=outer](4.98,0.84)(0.0,-0.84)
\pscircle[linewidth=0.02,dimen=outer](2.84,0.4){0.18}
\psdots[dotsize=0.12](2.84,-0.38)
\pscircle[linewidth=0.02,dimen=outer](2.84,-0.38){0.18}
\psline[linewidth=0.02cm](0.44,0.4)(0.44,0.42)
\psline[linewidth=0.02cm](0.48,0.42)(0.46,0.42)
\psline[linewidth=0.02cm](0.46,0.42)(0.44,0.4)
\psdots[dotsize=0.12](3.66,0.4)
\pscircle[linewidth=0.02,dimen=outer](3.66,0.4){0.18}
\psdots[dotsize=0.12](3.66,-0.38)
\pscircle[linewidth=0.02,dimen=outer](3.66,-0.38){0.18}
\psdots[dotsize=0.12](4.46,0.4)
\pscircle[linewidth=0.02,dimen=outer](4.46,0.4){0.18}
\psdots[dotsize=0.12](4.46,-0.38)
\pscircle[linewidth=0.02,dimen=outer](4.46,-0.38){0.18}
\psline[linewidth=0.02cm](2.82,0.42)(4.44,-0.38)
\psline[linewidth=0.02cm](3.66,0.42)(2.86,-0.36)
\psline[linewidth=0.02cm](4.46,0.4)(3.68,-0.36)
\end{pspicture}
\end{center}
\caption{a chain with three pairs of corresponding singletons}\label{fig3}
\end{figure}

\begin{figure}[h]
\begin{center}
\psset{unit=1cm}
\begin{pspicture}(0,-0.84)(1.8,0.84)
\psdots[dotsize=0.12](0.46,0.4)
\psdots[dotsize=0.12](1.24,0.4)
\psframe[linewidth=0.02,dimen=outer](1.8,0.84)(0.0,-0.84)
\psline[linewidth=0.02cm](0.44,0.4)(0.44,0.42)
\psline[linewidth=0.02cm](0.48,0.42)(0.46,0.42)
\psline[linewidth=0.02cm](0.46,0.42)(0.44,0.4)
\psellipse[linewidth=0.02,dimen=outer](0.87,0.41)(0.73,0.25)
\psdots[dotsize=0.12](0.46,-0.38)
\psdots[dotsize=0.12](1.24,-0.38)
\psline[linewidth=0.02cm](0.44,-0.38)(0.44,-0.36)
\psline[linewidth=0.02cm](0.48,-0.36)(0.46,-0.36)
\psline[linewidth=0.02cm](0.46,-0.36)(0.44,-0.38)
\psellipse[linewidth=0.02,dimen=outer](0.87,-0.37)(0.73,0.25)
\psline[linewidth=0.02cm](0.46,0.4)(0.46,-0.38)
\end{pspicture}
\end{center}
\caption{a chain with two corresponding blocks}\label{fig4}
\end{figure}
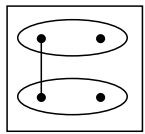
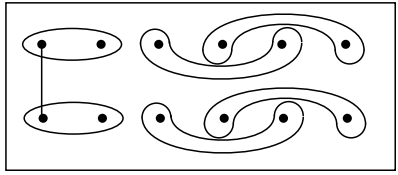
\begin{figure}[h]
\begin{center}
\psset{unit=1cm}
\begin{pspicture}(0,-1.12)(5.14,1.13)
\psdots[dotsize=0.12](3.64,0.56)
\psframe[linewidth=0.02,dimen=outer](5.14,1.12)(0.0,-1.12)
\psdots[dotsize=0.12](2.02,0.56)
\psdots[dotsize=0.12](1.26,0.56)
\psdots[dotsize=0.12](0.48,0.56)
\psellipse[linewidth=0.02,dimen=outer](0.88,0.56)(0.66,0.22)
\psbezier[linewidth=0.02](1.78,0.58)(1.78,-0.06)(3.9,-0.06)(3.9,0.58)
\psbezier[linewidth=0.02](2.16,0.6)(2.16,0.16)(3.52,0.16)(3.52,0.6)
\psarc[linewidth=0.02](1.97,0.57){0.19}{-0.0}{180.0}
\psarc[linewidth=0.02](3.71,0.59){0.19}{-0.0}{180.0}
\psdots[dotsize=0.12,dotangle=-180.0](2.86,0.56)
\psdots[dotsize=0.12,dotangle=-180.0](4.48,0.56)
\psbezier[linewidth=0.02](4.72,0.48)(4.72,1.12)(2.6,1.12)(2.6,0.48)
\psbezier[linewidth=0.02](4.34,0.48)(4.34,0.94)(2.98,0.94)(2.98,0.48)
\rput{-180.0}(9.06,0.98){\psarc[linewidth=0.02](4.53,0.49){0.19}{-0.0}{180.0}}
\rput{-180.0}(5.58,0.98){\psarc[linewidth=0.02](2.79,0.49){0.19}{-0.0}{180.0}}
\psdots[dotsize=0.12](3.66,-0.42)
\psdots[dotsize=0.12](2.04,-0.42)
\psdots[dotsize=0.12](1.28,-0.42)
\psdots[dotsize=0.12](0.5,-0.42)
\psellipse[linewidth=0.02,dimen=outer](0.9,-0.42)(0.66,0.22)
\psbezier[linewidth=0.02](1.8,-0.4)(1.8,-1.04)(3.92,-1.04)(3.92,-0.4)
\psbezier[linewidth=0.02](2.18,-0.38)(2.18,-0.82)(3.54,-0.82)(3.54,-0.38)
\psarc[linewidth=0.02](1.99,-0.41){0.19}{-0.0}{180.0}
\psarc[linewidth=0.02](3.73,-0.39){0.19}{-0.0}{180.0}
\psdots[dotsize=0.12,dotangle=-180.0](2.88,-0.42)
\psdots[dotsize=0.12,dotangle=-180.0](4.5,-0.42)
\psbezier[linewidth=0.02](4.74,-0.5)(4.74,0.14)(2.62,0.14)(2.62,-0.5)
\psbezier[linewidth=0.02](4.36,-0.5)(4.36,-0.04)(3.0,-0.04)(3.0,-0.5)
\rput{-180.0}(9.1,-0.98){\psarc[linewidth=0.02](4.55,-0.49){0.19}{-0.0}{180.0}}
\rput{-180.0}(5.62,-0.98){\psarc[linewidth=0.02](2.81,-0.49){0.19}{-0.0}{180.0}}
\psline[linewidth=0.02cm](0.44,0.56)(0.48,0.56)
\psline[linewidth=0.02cm](0.48,0.56)(0.48,-0.38)
\end{pspicture}
\end{center}
\caption{a class with three degrees of freedom}\label{fig5}
\end{figure}

\medskip\medskip
\noindent
{\bf 7.} Fix $k,N\geq 2$ and let $C\subset C(2k,N)$ be a generic subset of $C(2k,N)$. Let $q=1,...,2k$
be an integer. We say that $C$ has {\sl at most $q$ degrees of freedom} (or, equivalently, that $C$ {\sl has at
most $q$ free indices}) if there exists two subsets $D,E\subset [k]$ such that $|D|\geq 1$ and
the following two properties are verified: (i) $|D|+|E| \leq q $, and (ii) for every\footnote{As indicated by our notation, we regard $x_D$
and $y_E$ as vectors, respectively in $[N]^{|D|}$ and $[N]^{|E|}$, by endowing $D$ and $E$ with the natural
ordering induced by the ordering on $[k]$.} $x_D=\{x_a : a\in D\}\in
[N]^{|D|}$ and every $y_E=\{y_b : b\in E\}\in [N]^{|E|}$, there exists {\sl at most} one chain $c$ as in (\ref{chain}) such that $i_a = x_a$ for every $a\in D$ and $j_b = y_b$ for every $b\in E$. Note that our definition contemplates the possibility that $E = \emptyset$, and in this case the role of $y_E = \emptyset$ is immaterial. In other words, the class $C$ has at most $q$ degrees of freedom if every $c\in C$ is completely determined by those $i_a$ in the upper sub-chain such that $a\in D$ and those $j_b$ in the lower sub-chain such that $b\in E$. For instance, it is easily seen the class $C(2k,N)$
has (exactly) $2k$ degrees of freedom.
Another example is the diagram in Fig. \ref{fig5}, which corresponds
to the case $k=6$, $\pi  = \{\{1,2\},\{3,5\},\{4,6\}\}$ and $u=v=1$. One sees that, for every $N$,
specifying $i_1$, $i_4$ and $j_4$ completely identifies a chain inside the class
$C^{1,1}_\pi(12,N)$, which has therefore three degrees of
freedom.\footnote{Indeed, one has necessarily that $i_1=i_2=i_3=i_5=j_1=j_2=j_3=j_5$, $i_4=i_6$
and $j_4=j_6$.
}

\medskip

The proof of the two (useful) results contained in the next statement is elementary and omitted.

\begin{lemma} \label{L : easylemma} Fix $k,N\geq 2$.
\begin{itemize}
\item[{\rm (1)}] Let $q=1,...,2k$. Assume that a generic class $C\subset C(2k,N)$ has at most $q$ degrees of freedom. Then, $|C|\leq N^q$.
\item[{\rm (2)}] Let $\hat{1} = \{[k]\}$ be the one-block partition of $[k]$. Then, the class $C_{\hat{1} }(2k,N)$ contains only ``constant'' chains of the type (\ref{chain}) such that $(i_1,i_{2})=(i_a,i_{a+1})=(j_x,j_{x+1})$, for every $a=2,...,k$ and every $x=1,...,k$. It follows that $|C_{\hat{1}}(2k,N)| = N$.
\end{itemize}
\end{lemma}

Lemma \ref{L : easylemma} will be used in the subsequent section.

\subsection{Combinatorial upper bounds}\label{SS : CombBounds}

We keep the notation introduced in the previous section. The following statement, which is the key element for proving Proposition \ref{remainder-prop}, contains the main combinatorial estimate of the paper.

\begin{prop}\label{P : combUpper}
Fix $k,N\geq 2$, and let $\pi=\{B_1,...,B_r\}\in\mathcal{P}(k)$ be a partition containing at least one block of cardinality $\geq 2$. Let the class $C_\pi(2k,N)$ be defined as at Point {\bf 4}. of the previous section. Then, there exists a finite constant $\Theta(k,\pi)\geq 0$, depending only on $k$ and $\pi$ (and not on $N$), such that
\begin{equation}\label{CombFund}
|C_\pi(2k,N)| \leq \Theta(k,\pi) \times N^{k-1}.
\end{equation}
\end{prop}

\noindent{\it Proof}: We shall consider separately the two cases
\begin{itemize}
\item[{ A.}] For every $v=1,...,r$, $|B_v|\geq 2$.
\item[{ B.}] The partition $\pi$ contains at least one singleton.
\end{itemize}
 \underline{Case A.} When $k=2,3$, the only partition meeting the needed requirements is $\hat{1}$. According to Lemma \ref{L : easylemma}-(2), $|C_{\hat{1}}(2k,N) |=N$, so that the claim is proved, and we shall henceforth assume that $k\geq 4$. Start by observing that $r\leq k/2$. Moreover, the class $C_\pi(2k,N)$ contains only chains such that at least one block in the upper sub-chain corresponds to a block in the lower sub-chain, which yields in turn that
\begin{equation*}
C_\pi(2k,N) = \bigcup_{u,v=1}^r C^{u,v}_\pi(2k,N),
\end{equation*}
where we adopted the notation introduced at Point {\bf 5.} of Section \ref{SS : CombDefs}. This implies the crude estimate
\begin{equation}\label{easy}
|C_\pi(2k,N)|\leq \sum_{u,v=1}^r |C^{u,v}_\pi(2k,N)|.
\end{equation}
According to Lemma \ref{L : easylemma}-(1), it is now sufficient to prove that each class $C^{u,v}_\pi(2k,N)$ has at most $2r-1$ degrees of freedom: indeed, (\ref{easy}) together with the fact that $2r-1 \leq k-1$ would imply relation (\ref{CombFund}), with $\Theta(k,\pi) = r^2 \leq k^2/4$. Fix $u,v\in\{1,...,r\}$. To prove that $C^{u,v}_\pi(2k,N)$ has at most $2r-1$ degrees of freedom, we shall build two sets $D,E\subset [k]$ as follows. For every $s=1,...,r$, choose an element of the block $B_s$, and denote this element by $a_s$. Then, define
\[
D = \{a_s : s = 1,...,r\}, \,\,\, E=D \backslash \{a_v\},
\]
where `$ \backslash $' denotes the difference between sets. We now claim that, for every $x_D=\{x_a : a\in D\}\in [N]^{|D|}$ and every $y_E=\{y_b : b\in E\}\in [N]^{|E|}$, there exists at most one chain $c \in C^{u,v}_\pi(2k,N)$ as in (\ref{chain}) such that $i_a = x_a$ for every $a\in D$ and $j_b = y_b$ for every $b\in E$. To prove this fact, suppose that such a chain $c$ exists, and assume that there exists another chain
\[
c' = (i'_1,i'_2)(i'_2,i'_3)...(i'_k,i'_1)(j'_1,j'_2)(j'_2,j'_3)...(j'_k,j'_1)
\]
verifying this property and such that $c' \in C^{u,v}_\pi(2k,N)$.
The following hold:
(a) for every $s=1,...,r$ and every $a\in B_s$, one has that $i'_a = x_{a_s} = i_{a_s} = i_a$, (b) for every $s \neq v$ and every $a\in B_s$, $j'_a = y_{a_s} = j_{a_s} = j_a$ and (c) for $s=v$ and every $a\in B_v$,  \[j'_a = j'_{a_v} = i'_{a_u} = x_{a_u}=i_{a_u} = j_{a_v}=j_a.\]
As a consequence, $c'=c$. Since $|D|+|E| = 2r-1$, this concludes the proof of Proposition \ref{P : combUpper} in the Case A.

\bigskip

\noindent\underline{Case B}. We shall denote by $S$ the collection of the singleton(s) of $\pi$, that is
the subset of $[k]$ composed of those indices $a$ such that
$\{a\}\in\pi$. Note that $|S|>0$ by assumption. We also write $P$ for the collection of the indices $u\in [r]$ such that $|B_u|\geq 2$. Note that $P$ is a subset of $[r]$, whereas $S\subset [k]$. Note also that the set $[r]\backslash P$ is the collection of all those $v\in [r]$ such that $B_v$ is a singleton. Clearly, \[|P| = r-|S| \leq \frac{k-|S|}{2}.\] By exploiting the cyclic nature of sub-chains, we can always assume, without loss of generality, that $S$ contains the singleton $\{k\}$. Since $P$ is not empty, this entails that there exists at least one singleton of $\pi$ that is adjacent from the right to a block of cardinality at least two. Formally, this means that there exists $s^*\in S$ and $u^*\in P$ such that $s^*-1 \in B_{u^*}$. We shall distinguish two cases
\begin{itemize}
\item[{\rm B1.}] The block $B_{u^*}$ contains two consecutive integers.
\item[{\rm B2.}] The block $B_{u^*}$ does not contain two consecutive integers.
\end{itemize}

\begin{figure}[h]
\begin{center}
\psset{unit=1cm}
\begin{pspicture}(0,-1.05)(7.4,1.22)
\psframe[linewidth=0.02,dimen=outer](7.4,1.22)(0.0,-1.02)
\psdots[dotsize=0.12](2.88,0.66)
\psellipse[linewidth=0.02,dimen=outer](2.88,0.66)(0.2,0.2)
\psdots[dotsize=0.12](6.08,0.66)
\psellipse[linewidth=0.02,dimen=outer](6.08,0.66)(0.2,0.2)
\psdots[dotsize=0.12](2.08,0.66)
\psdots[dotsize=0.12](4.5,0.66)
\psdots[dotsize=0.12](3.68,0.66)
\psellipse[linewidth=0.02,dimen=outer](3.68,0.66)(0.2,0.2)
\psdots[dotsize=0.12](5.28,0.66)
\psbezier[linewidth=0.02](1.74,0.68)(1.74,-0.04)(4.82,-0.04)(4.82,0.68)
\psbezier[linewidth=0.02](2.26,0.7)(2.26,0.16)(4.34,0.16)(4.34,0.7)
\psarc[linewidth=0.02](2.0,0.68){0.26}{-0.0}{180.0}
\psline[linewidth=0.02cm](3.66,0.68)(3.68,0.66)
\psdots[dotsize=0.12](6.88,0.66)
\psellipse[linewidth=0.02,dimen=outer](6.88,0.66)(0.2,0.2)
\psarc[linewidth=0.02](4.91,0.43){0.63}{21.037512}{154.65382}
\psarc[linewidth=0.02](4.94,0.6){0.14}{26.565052}{151.69925}
\psarc[linewidth=0.02](5.28,0.8){0.26}{-149.03624}{-29.623749}
\psdots[dotsize=0.12](2.9,-0.34)
\psellipse[linewidth=0.02,dimen=outer](2.9,-0.34)(0.2,0.2)
\psdots[dotsize=0.12](6.1,-0.34)
\psellipse[linewidth=0.02,dimen=outer](6.1,-0.34)(0.2,0.2)
\psdots[dotsize=0.12](2.1,-0.34)
\psdots[dotsize=0.12](4.52,-0.34)
\psdots[dotsize=0.12](3.7,-0.34)
\psellipse[linewidth=0.02,dimen=outer](3.7,-0.34)(0.2,0.2)
\psdots[dotsize=0.12](5.3,-0.34)
\psbezier[linewidth=0.02](1.76,-0.32)(1.76,-1.04)(4.84,-1.04)(4.84,-0.32)
\psbezier[linewidth=0.02](2.28,-0.3)(2.28,-0.84)(4.36,-0.84)(4.36,-0.3)
\psarc[linewidth=0.02](2.02,-0.32){0.26}{-0.0}{180.0}
\psline[linewidth=0.02cm](3.68,-0.32)(3.7,-0.34)
\psdots[dotsize=0.12](6.9,-0.34)
\psellipse[linewidth=0.02,dimen=outer](6.9,-0.34)(0.2,0.2)
\psarc[linewidth=0.02](4.94,-0.58){0.64}{23.198591}{154.65382}
\psarc[linewidth=0.02](4.96,-0.4){0.14}{29.291363}{151.69925}
\psarc[linewidth=0.02](5.3,-0.2){0.26}{-149.03624}{-29.623749}
\psdots[dotsize=0.12](0.48,0.66)
\psellipse[linewidth=0.02,dimen=outer](0.48,0.66)(0.2,0.2)
\psdots[dotsize=0.12](1.28,0.66)
\psellipse[linewidth=0.02,dimen=outer](1.28,0.66)(0.2,0.2)
\psdots[dotsize=0.12](0.48,-0.34)
\psellipse[linewidth=0.02,dimen=outer](0.48,-0.34)(0.2,0.2)
\psdots[dotsize=0.12](1.28,-0.34)
\psellipse[linewidth=0.02,dimen=outer](1.28,-0.34)(0.2,0.2)
\end{pspicture}
\end{center}
\caption{a singleton is adjacent to a 3-block with two consecutive elements}\label{fig6}
\end{figure}
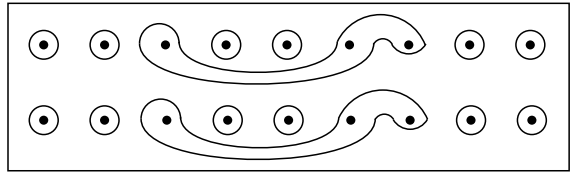

\noindent({\it Proof under B1}.) The situation of B1 is illustrated in Fig. \ref{fig6}, where $k=9$, \[\pi = \{B_1,...,B_7\}=\{\{1\},\{2\},\{3,6,7\},\{4\},\{5\},\{8\},\{9\}\},\] and one can take $s^*=8$, $u^*=3$, and the two consecutive integers in $B_{u^*}$ are 6 and 7.

\noindent Since each element of $C_\pi(2k,N)$ is such that every singleton in a given sub-chain corresponds to a block in the opposite sub-chain, we have that
\begin{equation}\label{UnionB1}
C_\pi(2k,N) = \bigcup_{R\in \mathcal{R}} C^R_\pi(2k,N),
\end{equation}
where we adopted the same notation as at Point {\bf 5.} of Section \ref{SS : CombDefs}, and the union runs over the class $\mathcal{R}$ of all bijections $R : U\rightarrow V $ such that both $U$ and $V$ contain the set $[r]\backslash P$, and every pair $(u,R(u))$ is such that at least one of the two blocks $B_u$ and $B_{R(u)}$ is a singleton. This entails the estimate
\begin{equation}\label{cool}
|C_\pi(2k,N)| \leq \sum_{R\in \mathcal{R}} |C^R_\pi(2k,N)|.
\end{equation}
To conclude the proof, we shall show that every class $C^R_\pi(2k,N)$ appearing in (\ref{cool}) has at most $k-1$ degrees of freedom: indeed, this fact together with Lemma \ref{L : easylemma}-(1) yields the desired conclusion (\ref{CombFund}), with the constant $\Theta(k,\pi) = |\mathcal{R}|$ (note that the definition of $\mathcal{R}$ does not depend on $N$) . To prove that $C^R_\pi(2k,N)$ has at most $k-1$ degrees of freedom, we define two sets $D,E\subset [k]$ as follows. For every $s=1,...,r$, choose an element of the block $B_s$, and denote this element by $a_s$. Then, define \[D = \{a_s : s = 1,...,r\}, \,\,\, E=D \backslash \left\{ \{a_{u^*}\}\cup \{a_s : s\in [r]\backslash P \}\right\}.\] In other words, $E$ is obtained by subtracting from $D$ the singleton(s) and the representative element of the block $B_{u^*}$, that is, of the block adjacent to $\{s^*\}$. We now want to prove that, for every $x_D=\{x_a : a\in D\}\in [N]^{|D|}$ and every $y_E=\{y_b : b\in E\}\in [N]^{|E|}$, there is at most one chain $c \in C^R_\pi(2k,N)$ as in (\ref{chain}) such that $i_a = x_a$ for every $a\in D$ and $j_b = y_b$ for every $b\in E$. To show this, assume that such a chain $c$ exists, and suppose that there exists another chain
\[
c' = (i'_1,i'_2)(i'_2,i'_3)...(i'_k,i'_1)(j'_1,j'_2)(j'_2,j'_3)...(j'_k,j'_1)
\]
verifying this property and such that $c' \in C^R_\pi(2k,N)$ and $c' \neq c$. By construction of the sets $D$ and $E$, all the indices composing the upper chain are completely determined by the choice of $x_D$, whereas the choice of $y_E$ determines the indices $j_x$ such that either $x$ is a singleton or $x\in B_v$ for some block $B_v$ of cardinality $\geq 2$ and such that $v\neq u^*$. This entails in turn that, necessarily since $c'\neq c$, one has that $j'_x \neq j_x$ for every $x\in B_{u^*}$. This is absurd. Indeed, since $B_{u^*}$ contains two consecutive integers, one has that $j'_x = j'_{x+1} $ and $j_x = j_{x+1}$ for every $x\in B_{u^*}$; it follows that, since $\{s^*\}$ is adjacent from the right to $B_{u^*}$ and therefore $s^*-1 \in B_{u^*}$ ,
\[
j'_x = j'_{s^*-1} = j'_{s^*} = y_{s*} =j_{s^*}= j_{s^*-1}=j_x,
\]
which is indeed a contradiction. Since \[|D|+|E| =r+|P|-1\leq \frac{k-|S|}{2} +|S| + \frac{k-|S|}{2} -1 =k-1,\] the proof is concluded.

\bigskip

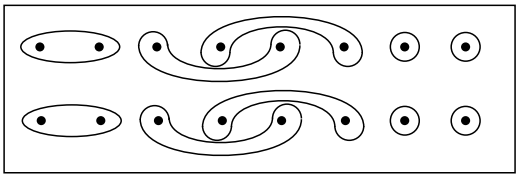
\begin{figure}[h]
\begin{center}
\psset{unit=1cm}
\begin{pspicture}(0,-1.12)(6.74,1.13)
\psdots[dotsize=0.12](3.64,0.56)
\psframe[linewidth=0.02,dimen=outer](6.74,1.12)(0.0,-1.12)
\psdots[dotsize=0.12](2.02,0.56)
\psdots[dotsize=0.12](5.28,0.56)
\psellipse[linewidth=0.02,dimen=outer](5.28,0.56)(0.2,0.2)
\psdots[dotsize=0.12](1.26,0.56)
\psdots[dotsize=0.12](6.08,0.56)
\psellipse[linewidth=0.02,dimen=outer](6.08,0.56)(0.2,0.2)
\psdots[dotsize=0.12](0.48,0.56)
\psellipse[linewidth=0.02,dimen=outer](0.88,0.56)(0.66,0.22)
\psbezier[linewidth=0.02](1.78,0.58)(1.78,-0.06)(3.9,-0.06)(3.9,0.58)
\psbezier[linewidth=0.02](2.16,0.6)(2.16,0.16)(3.52,0.16)(3.52,0.6)
\psarc[linewidth=0.02](1.97,0.57){0.19}{-0.0}{180.0}
\psarc[linewidth=0.02](3.71,0.59){0.19}{-0.0}{180.0}
\psdots[dotsize=0.12,dotangle=-180.0](2.86,0.56)
\psdots[dotsize=0.12,dotangle=-180.0](4.48,0.56)
\psbezier[linewidth=0.02](4.72,0.48)(4.72,1.12)(2.6,1.12)(2.6,0.48)
\psbezier[linewidth=0.02](4.34,0.48)(4.34,0.94)(2.98,0.94)(2.98,0.48)
\rput{-180.0}(9.06,0.98){\psarc[linewidth=0.02](4.53,0.49){0.19}{-0.0}{180.0}}
\rput{-180.0}(5.58,0.98){\psarc[linewidth=0.02](2.79,0.49){0.19}{-0.0}{180.0}}
\psdots[dotsize=0.12](3.66,-0.42)
\psdots[dotsize=0.12](2.04,-0.42)
\psdots[dotsize=0.12](5.28,-0.42)
\psellipse[linewidth=0.02,dimen=outer](5.28,-0.42)(0.2,0.2)
\psdots[dotsize=0.12](1.28,-0.42)
\psdots[dotsize=0.12](6.08,-0.42)
\psellipse[linewidth=0.02,dimen=outer](6.08,-0.42)(0.2,0.2)
\psdots[dotsize=0.12](0.5,-0.42)
\psellipse[linewidth=0.02,dimen=outer](0.9,-0.42)(0.66,0.22)
\psbezier[linewidth=0.02](1.8,-0.4)(1.8,-1.04)(3.92,-1.04)(3.92,-0.4)
\psbezier[linewidth=0.02](2.18,-0.38)(2.18,-0.82)(3.54,-0.82)(3.54,-0.38)
\psarc[linewidth=0.02](1.99,-0.41){0.19}{-0.0}{180.0}
\psarc[linewidth=0.02](3.73,-0.39){0.19}{-0.0}{180.0}
\psdots[dotsize=0.12,dotangle=-180.0](2.88,-0.42)
\psdots[dotsize=0.12,dotangle=-180.0](4.5,-0.42)
\psbezier[linewidth=0.02](4.74,-0.5)(4.74,0.14)(2.62,0.14)(2.62,-0.5)
\psbezier[linewidth=0.02](4.36,-0.5)(4.36,-0.04)(3.0,-0.04)(3.0,-0.5)
\rput{-180.0}(9.1,-0.98){\psarc[linewidth=0.02](4.55,-0.49){0.19}{-0.0}{180.0}}
\rput{-180.0}(5.62,-0.98){\psarc[linewidth=0.02](2.81,-0.49){0.19}{-0.0}{180.0}}
\end{pspicture}
\end{center}
\caption{a singleton is adjacent to a 2-block with no consecutive elements}\label{fig8}
\end{figure}

\noindent({\it Proof under B2}.) Since $B_{u^*}$ does not contain two consecutive integers and $|B_{u^*}|\geq 2$, we deduce the existence of a block $B_{\overline{u}}\in\pi$, which is different from $B_{u^*}$ and $\{s^*\}$, enjoying the following ``interlacement property'': there exists an integer $a\in[k]$ such that $a+1<s^*-1$, $a\in B_{u^*}$ and $a+1\in B_{\overline{u}}$. The block $B_{\overline{u}}$ can be either a singleton or a block with two or more elements. This situation is illustrated in Fig. \ref{fig8}, corresponding to the case $k=8$ and $\pi = \{B_1,...,B_5\}= \{\{1,2\},\{3,5\},\{4,6\},\{7\},\{8\}\}$. Here, $s^* = 7$, $B_{u^*} =B_3=\{4,6\}$, $B_{\overline{u}}=B_2=\{3,5\}$ and $a=4$.

\noindent The crucial remark is now that, for a chain $c$ as in (\ref{chain}) with partition $\pi$, one has that $i_{s^*} = i_{a+1}$. Indeed, $a$ and $s^*-1$ both belong to $B_{u^*}$, and therefore $ (i_{s^*-1},i_{s^*})= (i_a,i_{a+1})$. Since $a+1\in B_{\overline{u}}$, this fact yields in particular that, $i_x=i_{s^*}$ for every $x\in B_{\overline{u}}$, that is, the left indices associated with $B_{\overline{u}}$ are completely determined by the choice of $i_{s^*}$. By the same argument, one shows that $j_{s^*} = j_{a+1}$.  The rest of the proof is similar to the case B1. First, we observe that the representation (\ref{UnionB1}), with $\mathcal{R}$ defined exactly as for B1, continues to be true, from which we deduce the estimate (\ref{cool}). It is now sufficient to show that each class $C_\pi^R(2k,N)$ has at most $k-1$ degrees of freedom. To do this, one chooses a representative element from each block $B_s \in \pi$, noted $a_s$, and then defines the sets
\[D = \{a_s : s = 1,...,r, \, s\neq \overline{u}\}, \,\,\, E=D \backslash \left\{ a_s : s\in [r]\backslash P \right\},\]
that is, $D$ is built by selecting one element from each block of $\pi$, except for $B_{\overline{u}}$, and $E$ is obtained by subtracting from $D$ all the remaining indices $a$ such that $\{a\}$ is a singleton of $\pi$. One has that
\begin{equation}\label{nokia}
|D|+|E| \leq k-1.
\end{equation}
Indeed, $|D|= r-1=|P|+|S|-1 \leq  \frac{k-|S|}{2} +|S|-1$, and then one has to consider two cases: either (a) $B_{\overline{u}}$ is a singleton, from which it follows that $|E|=|D|- (|S|-1)\leq \frac{k-|S|}{2}$, or (b) $B_{\overline{u}}$ is not a singleton, yielding $|E|=|D|- |S|\leq \frac{k-|S|}{2}-1$. In these two cases, (\ref{nokia}) is then in order. To conclude, it remains to show that, for every $x_D=\{x_a : a\in D\}\in [N]^{|D|}$ and every $y_E=\{y_b : b\in E\}\in [N]^{|E|}$, there is at most one chain $c \in C^R_\pi(2k,N)$ as in (\ref{chain}) such that $i_a = x_a$ for every $a\in D$ and $j_b = y_b$ for every $b\in E$. To see this, assume that such a chain $c$ exists, and observe that, due to the above considerations, the choice of $x_D$ completely determines the upper sub-chain of $c$, as well as those indices $j_x$ in the lower sub-chain such that $\{x\}$ is a singleton of $\pi$ or (whenever $B_{\overline{u}}$ is not a singleton) such that  $x\in B_{\overline{u}}$. Since the remaining left indices in the lower sub-chain of $c$ are determined by the choice of $y_E$, the claim is proved. In view of (\ref{nokia}), this shows that $C_\pi^R(2k,N)$ has at most $k-1$ free indices. This concludes the proof of Proposition \ref{P : combUpper}.

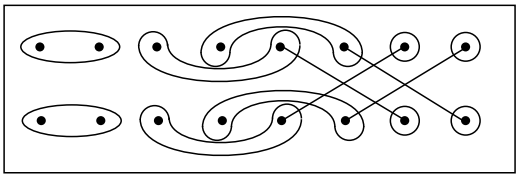
\begin{figure}[h]
\begin{center}
\psset{unit=1cm}
\begin{pspicture}(0,-1.12)(6.74,1.13)
\psdots[dotsize=0.12](3.64,0.56)
\psframe[linewidth=0.02,dimen=outer](6.74,1.12)(0.0,-1.12)
\psdots[dotsize=0.12](2.02,0.56)
\psdots[dotsize=0.12](5.28,0.56)
\psellipse[linewidth=0.02,dimen=outer](5.28,0.56)(0.2,0.2)
\psdots[dotsize=0.12](1.26,0.56)
\psdots[dotsize=0.12](6.08,0.56)
\psellipse[linewidth=0.02,dimen=outer](6.08,0.56)(0.2,0.2)
\psdots[dotsize=0.12](0.48,0.56)
\psellipse[linewidth=0.02,dimen=outer](0.88,0.56)(0.66,0.22)
\psbezier[linewidth=0.02](1.78,0.58)(1.78,-0.06)(3.9,-0.06)(3.9,0.58)
\psbezier[linewidth=0.02](2.16,0.6)(2.16,0.16)(3.52,0.16)(3.52,0.6)
\psarc[linewidth=0.02](1.97,0.57){0.19}{-0.0}{180.0}
\psarc[linewidth=0.02](3.71,0.59){0.19}{-0.0}{180.0}
\psdots[dotsize=0.12,dotangle=-180.0](2.86,0.56)
\psdots[dotsize=0.12,dotangle=-180.0](4.48,0.56)
\psbezier[linewidth=0.02](4.72,0.48)(4.72,1.12)(2.6,1.12)(2.6,0.48)
\psbezier[linewidth=0.02](4.34,0.48)(4.34,0.94)(2.98,0.94)(2.98,0.48)
\rput{-180.0}(9.06,0.98){\psarc[linewidth=0.02](4.53,0.49){0.19}{-0.0}{180.0}}
\rput{-180.0}(5.58,0.98){\psarc[linewidth=0.02](2.79,0.49){0.19}{-0.0}{180.0}}
\psdots[dotsize=0.12](3.66,-0.42)
\psdots[dotsize=0.12](2.04,-0.42)
\psdots[dotsize=0.12](5.28,-0.42)
\psellipse[linewidth=0.02,dimen=outer](5.28,-0.42)(0.2,0.2)
\psdots[dotsize=0.12](1.28,-0.42)
\psdots[dotsize=0.12](6.08,-0.42)
\psellipse[linewidth=0.02,dimen=outer](6.08,-0.42)(0.2,0.2)
\psdots[dotsize=0.12](0.5,-0.42)
\psellipse[linewidth=0.02,dimen=outer](0.9,-0.42)(0.66,0.22)
\psbezier[linewidth=0.02](1.8,-0.4)(1.8,-1.04)(3.92,-1.04)(3.92,-0.4)
\psbezier[linewidth=0.02](2.18,-0.38)(2.18,-0.82)(3.54,-0.82)(3.54,-0.38)
\psarc[linewidth=0.02](1.99,-0.41){0.19}{-0.0}{180.0}
\psarc[linewidth=0.02](3.73,-0.39){0.19}{-0.0}{180.0}
\psdots[dotsize=0.12,dotangle=-180.0](2.88,-0.42)
\psdots[dotsize=0.12,dotangle=-180.0](4.5,-0.42)
\psbezier[linewidth=0.02](4.74,-0.5)(4.74,0.14)(2.62,0.14)(2.62,-0.5)
\psbezier[linewidth=0.02](4.36,-0.5)(4.36,-0.04)(3.0,-0.04)(3.0,-0.5)
\rput{-180.0}(9.1,-0.98){\psarc[linewidth=0.02](4.55,-0.49){0.19}{-0.0}{180.0}}
\rput{-180.0}(5.62,-0.98){\psarc[linewidth=0.02](2.81,-0.49){0.19}{-0.0}{180.0}}
\psline[linewidth=0.02cm](3.64,0.58)(5.26,-0.4)
\psline[linewidth=0.02cm](4.48,0.6)(4.5,0.56)
\psline[linewidth=0.02cm](4.48,0.56)(6.08,-0.44)
\psline[linewidth=0.02cm](4.48,-0.4)(4.52,-0.42)
\psline[linewidth=0.02cm](4.5,-0.42)(6.1,0.56)
\psline[linewidth=0.02cm](3.64,-0.4)(3.66,-0.4)
\psline[linewidth=0.02cm](3.66,-0.42)(5.24,0.54)
\end{pspicture}
\end{center}
\caption{a class with three free indices}\label{fig9}
\end{figure}

\noindent As an illustration of the above arguments, one can consider the diagram in Fig. \ref{fig9}, that is constructed from the situation in Fig. \ref{fig8} by selecting $U=V=\{2,3,4,5\}$ and $R(2)=4$, $R(3)=5$, $R(4)=2$ and $R(5)=3$. In particular, it is easily seen that fixing $i_4$, $i_7$ and $i_8$ completely identifies a chain $c$ inside the class $C_\pi ^R (16,N)$, that has therefore three degrees of freedom.

\fin
\subsection{Proofs of Proposition \ref{remainder-prop} and Theorem \ref{Main}}\label{SS : Last Word}

{\it Proof of Proposition \ref{remainder-prop}}: We take up the notation introduced in Section \ref{SS : Sketch}. In view of Proposition \ref{P : combUpper}, in order to prove relation (\ref{gn}) (and therefore Proposition \ref{remainder-prop}), it is sufficient to show that, for every $\pi\in\mathcal{Q}(k)$, each pair $({\bf i},{\bf j})\in G_N(\pi)$ is such that the corresponding chain $(i_1,i_2)...(i_k,i_1)(j_1,j_2)...(j_k,j_1)$ is an element of $C_\pi (2k,N)$, from which one deduces $|G_N(\pi)|\leq |C_\pi (2k,N)|\leq\Theta(k,\pi)N^{k-1}$. To show the desired property, it is enough to prove that, for every pair $({\bf i},{\bf j})\in A_N(\pi)\times A_N(\pi)$ such that the chain $(i_1,i_2)...(i_k,i_1)(j_1,j_2)...(j_k,j_1)$ {\sl is not} in $C_\pi (2k,N)$, one has that $({\bf i},{\bf j})\not\in G_N(\pi)$.
By definition of $C_\pi(2k,N)$, we have to examine two cases.
Start by considering a partition $\pi\in\mathcal{Q}(k)$ not containing any singleton: if $({\bf i},{\bf j})\in A_N(\pi)\times A_N(\pi)$ is such that $(i_1,i_2)...(i_k,i_1)(j_1,j_2)...(j_k,j_1)\not\in C_\pi (2k,N)$, then the random variables $X_{i_ai_{a+1}}$ indexed by the upper sub-chain are independent of those indexed by the lower sub-chain, and consequently
\[
E(X_{i_1i_2}\ldots X_{i_ki_1}
X_{j_1j_2}\ldots X_{j_kj_1})
= E(X_{i_1i_2}\ldots X_{i_ki_1})
E(X_{j_1j_2}\ldots X_{j_kj_1}),
\]
yielding $({\bf i},{\bf j})\not\in G_N(\pi)$. On the other hand, if $\pi\in\mathcal{Q}(k)$ contains a singleton and if $({\bf i},{\bf j})$ is such that $(i_1,i_2)...(i_k,i_1)(j_1,j_2)...(j_k,j_1)\not\in C_\pi (2k,N)$, then there exists $a=1,...,k$ such that $X_{i_ai_{a+1}}$ or $X_{j_aj_{a+1}}$ is independent of all the other variables indexed by the elements of the chain. This gives
\[
E(X_{i_1i_2}\ldots X_{i_ki_1}
X_{j_1j_2}\ldots X_{j_kj_1})
= E(X_{i_1i_2}\ldots X_{i_ki_1})
E(X_{j_1j_2}\ldots X_{j_kj_1}) = 0,
\]
thus proving the required property $({\bf i},{\bf j})\not\in G_N(\pi)$. The proof is finished.
\fin

\noindent{\it Proof of Theorem \ref{Main}-(i)}: By virtue of the representation (\ref{representationINTRO})--(\ref{diagonals}) and of Proposition \ref{remainder-prop}, one sees that, for every $2\leq k_1<...<k_m$, the limit in distribution of the vector
\[
\Big({\rm Tr}(X_{N}),
{\rm Tr}(X^{k_1}_{N})-E\left[{\rm Tr}(X^{k_1}_{N})\right]\!
,\!...,
{\rm Tr}(X^{k_m}_{N})-E\left[{\rm Tr}(X^{k_m}_{N})\right]
\Big)
\]
coincides with the limit in distribution of

\[
\left(N^{-1/2}\sum_{i=1}^NX_{ii},\,\,\,
N^{-\frac{k_1}{2}} \!\!\!  \sum_{{\bf i}\in D_N^{(k_1)}} \!\!\! X_{i_1 i_2}X_{i_2 i_3}\cdot\cdot\cdot X_{i_{k_1} i_1}
,\,\ldots,
N^{-\frac{k_m}{2}} \!\!\!  \sum_{{\bf i}\in D_N^{(k_m)}} \!\!\!  X_{i_1 i_2}X_{i_2 i_3}\cdot\cdot\cdot X_{i_{k_m} i_1}
\right),
\]
so that the desired conclusion follows from Corollary \ref{P : megaP}.
\fin

\noindent{\it Proof of Theorem \ref{Main}-(ii)}:
For the simplicity of exposition, we assume that $k_1\geq 2$, the proof when
$k_1=1$ being completely similar and easier. We have, using the notation $D_N^{(k)}$ introduced
in the beginning of Section \ref{strategy} and using (\ref{diagonals}),
    \begin{eqnarray*}
&&  \Bigg|  E\left[\varphi\left(
\frac{{\rm Tr}(X^{k_1}_{N})-E[{\rm Tr}(X^{k_1}_{N})]}{\sqrt{{\rm Var}({\rm Tr}(X^{k_1}_{N}))}}
,\ldots,
\frac{{\rm Tr}(X^{k_m}_{N})-E[{\rm Tr}(X^{k_m}_{N})]}{\sqrt{{\rm Var}({\rm Tr}(X^{k_m}_{N}))}}
\right)\right]\\
&&\hskip7cm
- E\left[\varphi\left(\frac{Z_{k_1}}{\sqrt{k_1}},...,\frac{Z_{k_m}}{\sqrt{k_m}}\right)\right]\Bigg|\leq A_N
+B_N,
\end{eqnarray*}
where, by writing ${\rm Var}({\rm Tr}(X^{k_j}_{N})) = C_j(N)$,
    \begin{eqnarray*}
A_N&=&
  \Bigg|
E
\left[
\varphi
\left(
\frac{1}{C_1(N)^{1/2}N^{\frac{k_1}{2}}}\sum_{{\bf i}\in D_N^{(k_1)}}X_{i_1i_2}\ldots X_{i_{k_1}i_1}
,\ldots,
\right.\right.\\
&&\hskip2.0cm
\left.\left.\left.
\frac{1}{C_m(N)^{1/2}N^{\frac{k_m}{2}}}\sum_{{\bf i}\in D_N^{(k_m)}}X_{i_1i_2}\ldots X_{i_{k_m}i_1}
\right)\right]\right.
- E\left[\varphi\left(\frac{Z_{k_1}}{\sqrt{k_1}},...,\frac{Z_{k_m}}{\sqrt{k_m}}\right)\right]\Bigg|
\end{eqnarray*}
and
\begin{eqnarray*}
&&B_N = \\
&& \Bigg|
E
\left[
\varphi
\left(
\frac{1}{C_1(N)^{1/2}N^{\frac{k_1}{2}}}\sum_{{\bf i}\in D_N^{(k_1)}}X_{i_1i_2}\ldots X_{i_{k_1}i_1}
,\ldots,
\frac{1}{C_m(N)^{1/2}N^{\frac{k_m}{2}}}\sum_{{\bf i}\in D_N^{(k_m)}}X_{i_1i_2}\ldots X_{i_{k_m}i_1}
\right)\right]\\
&&\hskip2cm
-
E\left[\varphi\left(\frac{{\rm Tr}(X^{k_1}_{N})-E[{\rm Tr}(X^{k_1}_{N})]}{\sqrt{{\rm Var}({\rm Tr}(X^{k_1}_{N}))}}
,\ldots,
\frac{{\rm Tr}(X^{k_m}_{N})-E[{\rm Tr}(X^{k_m}_{N})]}{\sqrt{{\rm Var}({\rm Tr}(X^{k_m}_{N}))}}
\right)\right]\Bigg|.
\end{eqnarray*}
By combining Corollary \ref{cocominet} with the computations made in the proof of Proposition \ref{1dim},
we immediately get
that $A_N=O(N^{-1/4})$. For $B_N$, we can write
\begin{eqnarray*}
|B_N|&\leq&
K\|\varphi'\|_\infty
\sum_{j=1}^m
E\left[
N^{-\frac{k_j}{2}}\left|\sum_{{\bf i}\not\in D_N^{(k_j)}}\big(X_{i_1i_2}\ldots X_{i_{k_j}i_1}
-E[X_{i_1i_2}\ldots X_{i_{k_j}i_1}]\big)
\right|\right]\\
&\leq&
K\|\varphi'\|_\infty
\sum_{j=1}^m
\sqrt{{\rm Var}\left(
N^{-\frac{k_j}{2}}\left|\sum_{{\bf i}\not\in D_N^{(k_j)}}\big(X_{i_1i_2}\ldots X_{i_{k_j}i_1}
-E[X_{i_1i_2}\ldots X_{i_{k_j}i_1}]\big),
\right|\right)},
\end{eqnarray*}
for some constant $K$ not depending on $N$, so that $B_N=O(N^{-1/2})=O(N^{-1/4})$ by Proposition \ref{remainder-prop}.
\fin

\section{Almost sure central limit theorems (ASCLTs)}\label{S : ASCLT}
\subsection{Preliminaries: a result by Ibragimov and Lifshits}
For $x,y\in\R^m$ ($m\geq 1$ fixed),
we write $\langle x,y\rangle=x_1y_1+\ldots+x_my_m$ (resp. $|x|=\sqrt{\langle x,x\rangle}$)
to indicate the inner product of $x$ and $y$ (resp. the norm of $x$). The following result, due to Ibragimov and Lifshits, plays a crucial role in the proof of Theorem \ref{thm-ASCLT}.

\begin{thm}[See \cite{IL}]\label{thm-IL}
Let $G = \{G_n : n\geq 1\}$ be a sequence of $\R^m$-valued random variables converging in distribution
towards a random
variable $G_\infty$, and set
\begin{equation}\label{delta}
\Delta_N(G,t)=\frac1{\log N}\sum_{n=1}^N \frac{1}n \big(e^{i\langle t,G_n\rangle}-
E[e^{i\langle t,G_\infty\rangle}]\big),\quad t\in\R^m.
\end{equation}
If, for  all $r>0$,
\begin{equation}\label{cond-IL}
\sup_{|t|\leq r}\sum_{N=2}^\infty \frac{E\vert \Delta_N(G,t)\vert^2}{N\log N}<\infty,
\end{equation}
then, almost surely, for all continuous and bounded function
$\varphi:\R^m\to\R$, we have
\begin{equation}\label{cl-IL}
\frac{1}{\log N}\sum_{n=1}^{N} \frac{\varphi(G_n)}{n}
\longrightarrow E[\varphi(G_\infty)],\quad\mbox{as $N\to\infty$}.
\end{equation}
\end{thm}
\begin{rem}
{\rm
\begin{enumerate}
\item If $E|\Delta_N(G,t)|^2=O(1/\log N)$ uniformly in $t$
on bounded sets, then (\ref{cond-IL})
is automatically satisfied.
\item See \cite{BNT} for several applications of Theorem \ref{thm-IL} in the framework of ASCLTs on Wiener space.
\end{enumerate}
}
\end{rem}

The following useful result allows to deal with sequences of random variables having the form of a sum of two terms, one of which vanishes in the mean-square sense.

\begin{lemme}\label{reste}
Let $G = \{G_n:n\geq 1 \}$ be a sequence of $\R^m$-valued random variables converging in distribution
towards a random
variable $G_\infty$, and satisfying in addition (\ref{cond-IL}).
Let $R = \{R_n:n\geq 1 \}$ be a sequence of $\R^m$-valued random variables
converging in $L^2(\Omega)$ to $R_\infty=0$, and such that
\begin{equation}\label{cond-reste}
\sum_{N=2}^\infty \frac{1}{N\log^2 N}\sum_{n=1}^N \frac{1}n E|R_n|^2<\infty.
\end{equation}
Then
\[
\sup_{|t|\leq r}\sum_{N=2}^\infty \frac{
E\big|\Delta_N(G+R,t)\big|^2
}{N\log N}<\infty,
\]
where $G+R = \{G_n+R_n : n\geq 1\} $ and $\Delta_N(G+R,t)$ is defined according to (\ref{delta}).
\end{lemme}
\begin{rem}
{\rm
If $E|R_n|^2=O(n^{-a})$, for some $a>0$, then (\ref{cond-reste})
is automatically satisfied.
}
\end{rem}
{\it Proof of Lemma \ref{reste}}. Since
$\Delta_N(G+R,t)=\sum_{n=1}^N \frac{1}n \big(e^{i\langle t,G_n+R_n\rangle}-
E[e^{i\langle t,G_\infty\rangle}]\big)$, one has that
\[
\Delta_N(G+R,t)=
\Delta_N(G,t)+
\frac1{\log N}\sum_{n=1}^N \frac{1}n e^{i\langle t,G_n\rangle}\big(
e^{i\langle t,R_n\rangle}-1
\big),
\]
so that, by using $|x+y|^2\leq 2|x|^2+2|y|^2$, Jensen inequality and
$\sum_{n=1}^N\frac1n\sim \log N$ as $N\to\infty$,
there exists a constant $c>0$ (independent of $N$) such that, for all $N\geq 2$,
\[
E\vert \Delta_N(G+R,t)\vert^2\leq 2
E\vert \Delta_N(G,t)\vert^2
+\frac{c}{\log N}
\sum_{n=1}^N \frac{1}n
E\big|e^{i\langle t,R_n\rangle}-1\big|^2.
\]
Since $|e^{i\langle t,x\rangle}-1|\leq |t||x|$, we deduce
\[
E\vert \Delta_N(G+R,t)\vert^2
\leq
2E\vert \Delta_N(G,t)\vert^2
+\frac{c|t|^2}{\log N}
\sum_{n=1}^N \frac{1}n
E|R_n|^2.
\]
The desired conclusion follows.\fin

\subsection{Proof of Theorem \ref{thm-ASCLT}}

For the sake of brevity, we shall prove Theorem \ref{thm-ASCLT} only for powers $k_i$ {\sl strictly}
greater than one. The general case ($k_i\geq 1$) can be deduced from similar arguments.

\medskip

\noindent Throughout this section, we fix integers $m\geq 1$ and $k_m>\ldots>k_1\geq 2$.
For $N\geq 1$ and $k\geq 2$, we denote (as above) by $D_N^{(k)}$ the collection
of all vectors ${\bf i}=(i_1,\ldots,i_{k})\in\{1,\ldots,N\}^k$ such that all pairs $(i_a,i_{a+1})$,
$a=1,\ldots,k$, are different (with the convention that $i_{k+1}=i_1$), that is,
${\bf i}\in D_N^{(k)}$ if and only if $(i_a,i_{a+1})\neq (i_b,i_{b+1})$ for every $a\neq b$ and $1\leq i_a \leq N$ for every $a=1,...,k$.
Let
\[
J_N(k)=N^{-k/2}\sum_{{\bf i}\in D_N^{(k)}}X_{i_1i_2}X_{i_2i_3}\ldots X_{i_ki_1},
\quad\mbox{and}\quad L_N(k)=\frac{J_N(k)}{\sqrt{E[J_N(k)^2]}}.
\]
Observe that $E[J_N(k)]=E[L_N(k)]=0$ and ${\rm Var}[L_N(k)]=1$. The proof of Theorem
\ref{thm-ASCLT} is divided into several steps.\\
\\
{\it \underline{Step 1}: bounding $E[L_n(k)L_p(k)]$}.
Fix $k\geq 2$. We shall prove
that there exists a constant $C_k>0$ such that, for all $n,p\geq 1$,
\begin{equation}\label{step1}
E[L_n(k)L_p(k)]\leq C_k\sqrt{\frac{n\wedge p}{n\vee p}}.
\end{equation}
By symmetry, we assume without loss of generality that
$p\geq n$.
If ${\bf i}\in D_n^{(k)}$ and ${\bf j}\in D_p^{(k)}\setminus D_n^{(k)}$,
then
\[
E[X_{i_1i_2}\ldots X_{i_ki_1}X_{j_1j_2}\ldots X_{j_kj_1}] = 0;
\]
indeed, if ${\bf i}\in D_n^{(k)}$ and ${\bf j}\in D_p^{(k)}\setminus D_n^{(k)}$, then necessarily there exists $a=1,...,k$ such that $j_a >n$, and therefore the centered random variable $X_{j_aj_{a+1}}$ is independent of $X_{i_b i_{b+1}}$ for every $b=1,...,k$, and also (by the definition of $D_p^{(k)}$) independent of $X_{j_sj_{s+1}}$ for every $s\neq a$. It follows that
\begin{eqnarray*}
E[J_n(k)J_p(k)]&=&
\left(\frac{1}{np}\right)^{k/2}\sum_{{\bf i}\in D_n^{(k)}}
\sum_{{\bf j}\in D_p^{(k)}} E[X_{i_1i_2}\ldots X_{i_ki_1}X_{j_1j_2}\ldots X_{j_kj_1}]\\
&=&\left(\frac{1}{np}\right)^{k/2}\sum_{{\bf i}\in D_n^{(k)}}
\sum_{{\bf j}\in D_n^{(k)}} E[X_{i_1i_2}\ldots X_{i_ki_1}X_{j_1j_2}\ldots X_{j_kj_1}]\\
&=&\left(\frac{n}{p}\right)^{k/2}E\left[J_n(k)^2\right].
\end{eqnarray*}
Thus
$
E[L_n(k)L_p(k)]=\left(\frac{n}{p}\right)^{k/2}\sqrt{\frac{E[J_n(k)^2]}{E[J_p(k)^2]}}.
$
But we have
$E[J_n(k)^2]\to k$ as $n\to\infty$, see indeed (\ref{prop1}).
As a consequence, we immediately
get the existence of a constant $C_k$ such that (\ref{step1}) is in order.\\
\\
{\it \underline{Step 2}: showing that $\sup_{|t|\leq r}\sum_{N=2}^\infty \frac{E\vert \Delta_N(L,t)\vert^2}{N\log N}<\infty$}.
Fix $k\geq 2$. Let $f_{k,N}$ be as in (\ref{zidane}).
Set $g_{k,N}=\frac{1}{\sqrt{E[J_N(k)^2]}}f_{k,N}$. We obviously
have $L_N(k)=Q_k(g_{k,N},{\bf X})$.
Combining (\ref{prop1}) and (\ref{speed}),
we immediately get that, for all $r=1,\ldots,k-1$,
$\|g_{k,N}\star_r g_{k,N}\|_{2k-2r}=O(N^{-1/2})$.
From now on,  for simplicity write
$L_N=\big(L_N(k_1),\ldots,L_N(k_m)\big)$, $N\geq 1$, and $g(t)=e^{-|t|^2/2}$, $t\in\R^m$.
Corollary \ref{cocominet} yields that
\begin{equation}\label{nancy1}
\left|E[e^{i\langle t,L_N\rangle}]-g(t)\right|=O(N^{-1/4}).
\end{equation}
On the other hand, for all $r=1,\ldots,k-1$, we can write
\begin{eqnarray*}
&&\|(g_{k,N}-g_{k,M})\star_r (g_{k,N}-g_{k,M})\|_{2k-2r}\\
&=&
\|g_{k,N}\star_r g_{k,N}+g_{k,M}\star_r g_{k,M}-g_{k,N}\star_r g_{k,M}-g_{k,M}\star_r
g_{k,N}\|_{2k-2r}
\\
&\leq&
\|g_{k,N}\star_r g_{k,N}\|_{2k-2r}+\|g_{k,M}\star_r g_{k,M}\|_{2k-2r}+2
\|g_{k,N}\star_r
g_{k,M}\|_{2k-2r}.
\end{eqnarray*}
But
\begin{eqnarray*}
\|g_{k,N}\star_r
g_{k,M}\|_{2k-2r}&=&
\sqrt{\big\langle g_{k,N}\star_{k-r} g_{k,N}, g_{k,M}\star_{k-r}g_{k,M}\big\rangle_{2r}}\\
&\leq&
\sqrt{\| g_{k,N}\star_{k-r} g_{k,N}\|_{2r}}\sqrt{
\| g_{k,M}\star_{k-r}g_{k,M}\|_{2r}}
\\
&\leq&\frac12\big(
\| g_{k,N}\star_{k-r} g_{k,N}\|_{2r}+
\| g_{k,M}\star_{k-r}g_{k,M}\|_{2r}
\big).
\end{eqnarray*}
Consequently,
\[
\left\|\frac{g_{k,N}-g_{k,M}}{\sqrt{2}}\star_r \frac{g_{k,N}-g_{k,M}}{\sqrt{2}}\right\|_{2k-2r}
\leq
\|g_{k,N}\star_r g_{k,N}\|_{2k-2r}+\|g_{k,M}\star_r g_{k,M}\|_{2k-2r}
=O(N^{-1/2}),
\]
as $N\to\infty$, uniformly on $M\geq N$ and $r=1,...,k-1$, that is, there exists a constant $C_k>0$ (depending solely on $k$) such that, for every $N$,
\[
\sup_{1\leq r\leq k-1\, ; \,  M\geq N}\left\|\frac{g_{k,N}-g_{k,M}}{\sqrt{2}}\star_r \frac{g_{k,N}-g_{k,M}}{\sqrt{2}}\right\|_{2k-2r} \leq\,\, \frac{C_k}{N}\,.
\]
Since ${\rm Var}\left[\frac{L_N(k)-L_M(k)}{\sqrt{2}}\right]=1-E[L_N(k)L_M(k)]$
for all $k$, by using Corollary \ref{cocominet} with \[
Q^i_{M,N}({\bf X})
=\big(L_N(k_i)-L_M(k_i)\big)/\sqrt{2-2E[L_N(k_i)L_M(k_i)]}\]
and
\[
\varphi(z_1,\ldots,z_m)=\exp\left(
\sqrt{1-E[L_N(k_1)L_M(k_1)]}t_1z_1
+\ldots+
\sqrt{1-E[L_N(k_m)L_M(k_m)]}t_mz_m
\right),
\]
we get that

\begin{equation}\label{nancy2}
\left|E\left[e^{i\left\langle t,\frac{L_N-L_M}{\sqrt{2}}\right\rangle}\right]-
\exp\left(-\big(1-E[L_N(k_1)L_M(k_1)]\big)\frac{t_1^2}{2}-\ldots-\big(1-
E[L_N(k_m)L_M(k_m)]\big)\frac{t_m^2}2\right)
\right|
\end{equation}
is $O(N^{-1/4})$
as $N\to\infty$, uniformly on $M\geq N$.
On the other hand,
combining (\ref{step1}) with
$\big|e^{-x^2/2}-e^{-(1-\alpha)x^2/2}\big|\leq \alpha x^2/2$ for all $x\in\R$ and $\alpha\geq 0$,
we get
that there exists $C_r>0$ such that, for all $t\in\R^m$ with $|t|\leq r$,
\begin{eqnarray}\notag
&& \left|g(t)-
\exp\left(-\big(1-E[L_N(k_1)L_M(k_1)]\big)\frac{t_1^2}{2}-\ldots-\big(1-
E[L_N(k_m)L_M(k_m)]\big)\frac{t_m^2}2\right)\right|\\
&& \leq C_r\sqrt{\frac{N\wedge M}{N\vee M}}. \label{nancy3}
\end{eqnarray}
Define $\Delta_N(L,t)$ according to (\ref{delta}), with $L_\infty\sim \mathscr{N}_m(0,I_m)$.
For $|t|\leq r$, we have, due to (\ref{nancy1})-(\ref{nancy2})-(\ref{nancy3}):
\begin{eqnarray*}
&\!\!&E\vert \Delta_N(L,t)\vert^2 \\
&\!=\!&\frac{1}{\log^2 N}\sum_{n,p=1}^N \frac{1}{np}
E\left[\big(e^{i\langle t,L_n\rangle}-g(t)\big)\big(e^{-i\langle t,L_p\rangle}-g(t)\big)
\right]\\
&\!=\!&\frac{1}{\log^2 N} \sum_{n,p=1}^{N}\frac{1}{np}
\left[
\left(E\big(e^{i\langle t,L_n-L_p\rangle}\big)-g^2(t)\right)
- g(t)\left(
E\big(e^{i\langle t,L_n\rangle}\big)-g(t)
\right)\right.\\
&&\hskip9cm\left.
- g(t)\left(
E\big(e^{-i\langle t,L_p\rangle}\big)-g(t)
\right)
\right]\notag\\
&\!=\!&\frac{1}{\log^2 N} \sum_{n,p=1}^{N}\frac{1}{np}
\left[
\left(
E\big(
e^{
i\sqrt{2}\left\langle
t,\frac{L_n-L_p}{\sqrt{2}}
\right\rangle
}
\big) - g(\sqrt{2}\,t)\right)
- g(t)\left(
E\big(e^{i\langle t,L_n\rangle}\big)-g(t)
\right)\right.\\
&&\hskip9cm\left.
- g(t)\left(
E\big(e^{-i\langle t,L_p\rangle}\big)-g(t)
\right)
\right].\notag\\
&\leq&\frac{C_r}{\log^2 N} \sum_{n,p=1}^{N}\frac{1}{np}
\left(
\sqrt{\frac{n\wedge p}{n\vee p}}+\frac1{n^{1/4}}+\frac1{p^{1/4}}
\right).
\end{eqnarray*}
It is obvious that
\[
\frac{1}{\log^2 N}\sum_{n,p=1}^N\frac{1}{np}\left(\frac{1}{n^{1/4}}+\frac{1}{p^{1/4}}\right)
=O\left(\frac{1}{\log N}\right)\quad \mbox{as $N\to\infty$}.
\]
Moreover,
\begin{eqnarray*}
\frac1{\log^2 N}\sum_{n,p=1}^N \frac1{np}\sqrt{\frac{n\wedge p}{n\vee p}}&\leq&
\frac{2}\log^2 N\sum_{n=1}^N \frac{1}{n\sqrt{n}}\sum_{p=1}^n \frac{1}{\sqrt{p}}\\
&\leq&\frac{c}{\log^2 N}\sum_{n=1}^N \frac{1}{n} = O\left(\frac1{\log N}\right).
\end{eqnarray*}
Hence, $\sup_{|t|\leq r}E\vert \Delta_N(L,t)\vert^2 =O(1/\log N)$,
implying immediately
\begin{equation}\label{nancy4}
\sup_{|t|\leq r}\sum_{N=2}^\infty \frac{E\vert \Delta_N(L,t)\vert^2}{N\log N}<\infty.
\end{equation}
{\it \underline{Step 3}: using Lemma \ref{reste}.}
Set $T_N(k)=\frac{J_N(k)}{\sqrt{k}}$. Using (\ref{prop1})
and elementary calculations,
it is immediate that
\[
E\left|
T_N(k)-L_N(k)
\right|^2 =
\frac{
\left|
E\left[J_N(k)^2\right]-k
\right|^2
}{
k
\left(\sqrt{k}+\sqrt{
E\left[J_N(k)^2\right]
}
\right)^2}\leq\frac{C_k^2}{k^2\,N^2},
\]
so that $\widetilde{R}_N:=T_N-L_N$ verifies condition (\ref{cond-reste}) of Lemma \ref{reste}.
Since (\ref{nancy4}) is also in order, we deduce that
$\sup_{|t|\leq r}\sum_{N=2}^\infty \frac{E\vert \Delta_N(T,t)\vert^2}{N\log N}<\infty$,
which in turns implies
\begin{equation}\label{nancy5}
\sup_{|t|\leq r}\sum_{N=2}^\infty \frac{E\vert \Delta_N(J,t)\vert^2}{N\log N}<\infty.
\end{equation}
{\it \underline{Step 4}: using Lemma \ref{reste} once again.}
For any $k\geq 2$, set $S_N(k)=
{\rm Tr}(A_N^k) - E\left[
{\rm Tr}(A_N^k)
\right]$.
We have
\[
S_N(k)=
N^{-\frac{k}{2}}  \sum_{i_1,...,i_k=1}^N
\big(X_{i_1 i_2}X_{i_2 i_3}\cdot\cdot\cdot X_{i_k i_1}
 - E[
X_{i_1 i_2}X_{i_2 i_3}\cdot\cdot\cdot X_{i_k i_1}
]\big)=
J_N(k)+R_N(k),
\]
with
\[
R_N(k)=N^{-\frac{k}{2}}  \sum_{{\bf i}\not\in D_N^{(k)}}
\big(X_{i_1 i_2}X_{i_2 i_3}\cdot\cdot\cdot X_{i_k i_1}
 - E[
X_{i_1 i_2}X_{i_2 i_3}\cdot\cdot\cdot X_{i_k i_1}
]\big).
\]
For all $k\geq 2$, we have
$E|R_N(k)|^2=O(1/N)$, see (\ref{Graal}).
Hence, Lemma \ref{reste} together with (\ref{nancy5}) imply that
$\sup_{|t|\leq r}\sum_{N=2}^\infty \frac{E\vert \Delta_N(S,t)\vert^2}{N\log N}<\infty$.
To finish the proof of Theorem \ref{thm-ASCLT}, it suffices to apply
Theorem \ref{thm-IL}.
\fin

\bigskip

 \bigskip

\noindent {\bf Acknowledgments.} We thank Paul Bourgade, Mikhail Lifshits, Gesine
Reinert and Brian Rider for helpful discussions.


\begin{thebibliography}{99}

\bibitem{AirMalVie} H. Airault, P. Malliavin and F. Viens (2009).
Stokes formula on the Wiener space and $n$-dimensional Nourdin-Peccati analysis.
{\it J. Funct. Anal.} {\bf 248}, no. 5, 1763-1783.

\bibitem{AGZ} G. Anderson, A. Guionnet and O. Zeitouni (2009). An introduction to random matrices.
Cambridge Studies in Advanced Mathematics {\bf 118}.

\bibitem{AZ} G. Anderson and O. Zeitouni (2006). A CLT for a band matrix model.
{\it Probab. Theory Related Fields} {\bf 134}(2), 283-338.

\bibitem{BaiYin} Z.D. Bai and Y.Q. Yin (1986). Limiting behaviour of the norm products of random matrices and two problems of Geman-Hwang. {\it Probab. Theory Rel. Fields} {\bf 73}, 555--569.

\bibitem{BNT} B. Bercu, I. Nourdin and M.S. Taqqu (2009).
A multiple stochastic integral criterion for almost sure limit theorems.
In revision for: {\it Stoch. Proc. Appl.}

\bibitem{BerkesCsaki}
I. Berkes and E. Cs\'aki (2001).
A universal result in almost sure central limit theory.
{\it Stoch. Proc. Appl.} {\bf 94}, no. 1, 105-134.

\bibitem{Brosamler}
G. A. Brosamler (1988).
An almost everywhere central limit theorem.
{\it Math. Proc. Cambridge Philos. Soc.} {\bf 104}, no. 3, 561-574.

\bibitem{Chatterjee_ptrf} S. Chatterjee (2009). Fluctuation of eigenvalues and second order Poincaré inequalities.
 \textit{Probab. Theory Rel. Fields} {\bf 143}, 1--40.

\bibitem{chen-shao}
\rm L.H.Y. Chen and Q.-M. Shao (2005).
\rm Stein's method for normal approximation.
In: {\it An Introduction to Stein's Method} (A.D. Barbour and L.H.Y. Chen, eds), Lecture Notes Series No.{\bf  4},
Institute for Mathematical Sciences, National University of Singapore, Singapore University Press and World Scientific 2005, 1--59.

\bibitem{CostLeb} A. Costin and J.L. Lebowitz (1995). Gaussian fluctuations in random matrices. {\it Physical Review
Letters} {\bf 75}, 69--72.

\bibitem{DE} P. Diaconis and S.E. Evans (2001). Linear functionals of eigenvalues of random matrices.
{\it Trans. Amer. Math. Soc.} {\bf 353}, no. 7, 2615-2633.

\bibitem{DS} P. Diaconis and M. Shahshahani (1994). On the eigenvalues of random matrices.
{\it J. Appl. Probab.} {\bf 31}, 49-62.


\bibitem{forrester1999} P.J. Forrester (1999). Fluctuation formula for complex random matrices. {\it J. Phys.
A} {\bf 32}, 159--163.

\bibitem{Geman80} S. Geman (1980). A limit theorem for the norm of random matrices. {\it Ann. Probab.}
{\bf 8} 252--261.

\bibitem{Geman1986} S. Geman (1986). The spectral radius of large random matrices. {\it  Ann. Probab.} {\bf 14}
1318--1328.

\bibitem{Ginibre1965} J. Ginibre (1965). Statistical ensembles of complex, quaternion and real matrices. {\it J. Math. Phys.} {\bf 6}, 440--449

\bibitem{GuionnetBook} A. Guionnet (2008). {\it Large random matrices: lectures on macroscopic asymptotics}. LNM {\bf 1957}. Springer.

\bibitem{hormann}
S. H\"ormann (2007):
Critical Behavior in Almost Sure Central Limit Theory.
{\it J. Theoret. Probab.} {\bf 20}, 613-636.

\bibitem{IL}
I. A. Ibragimov and M. A. Lifshits (2000).
On limit theorems of ``almost sure'' type.
{\it Theory Probab. Appl.} {\bf 44}, no. 2, 254-272.

\bibitem{Janson} S. Janson (1997). \textit{Gaussian Hilbert Spaces.} Cambridge University Press, Cambridge.

\bibitem{LaceyPhillip}
M.T. Lacey and W. Philipp (1990).
A note on the almost sure central limit theorem.
{\it Statist. Probab. Letters} {\bf 9}, 201-205.

\bibitem{Levy}
P. L\'evy (1937).
{\it Th\'eorie de l'addition des variables al\'eatoires}.
Gauthiers-Villars.


\bibitem{MOO} E. Mossel, R. O'Donnell and K. Oleszkiewicz (2010).
Noise stability of functions with low influences: invariance and optimality.
{\it Ann. Math.} {\bf 171}, no. 1, 295-341.

\bibitem{NP-PTRF}
\rm I. Nourdin and G. Peccati (2009).
\rm Stein's method on Wiener chaos.
{\it Probab. Theory Rel. Fields} {\bf 145}, no. 1, 75--118.


\bibitem{Noupecrei3}
\rm I. Nourdin, G. Peccati and G. Reinert (2009).
\rm Invariance principles for homogeneous sums: universality of Gaussian Wiener chaos.
\rm To appear in: {\it Ann. Probab.}


\bibitem{nualartbook}
\rm D. Nualart (2006).
\it The Malliavin calculus and related topics.
\rm Springer Verlag, Berlin, Second edition, 2006.

\bibitem{nunugio}
\rm D. Nualart and G. Peccati (2005).
\rm Central limit theorems for sequences of multiple stochastic integrals.
{\it Ann. Probab.} {\bf 33} (1), 177--193.

\bibitem{PTu04} G.\ Peccati and C.A. Tudor (2005). Gaussian limits for
vector-valued multiple stochastic integrals. \textit{S\'{e}minaire de
Probabilit\'{e}s XXXVIII}, LNM \textbf{1857}. Springer-Verlag, Berlin
Heidelberg New York, pp. 247--262.

\bibitem{riderptrf2004} B. Rider (2004).
Deviations from the circular law. {\it Probab. Theory Related Fields}
{\bf 130}, 337--367

\bibitem{RiderSilvAop2006} B. Rider and J. Silverstein (1986).
Gaussian fluctuations for non-Hermitian random matrix ensembles. {\it Ann. Probab.} {\bf 34}, no. 6, 2118--2143

\bibitem{RiderVirag} B. Rider and B. Vir\'{a}g (2007). The noise in the circular law and the
Gaussian free field. {\it Int. Math. Res. Not.} 2, Art. ID rnm006.

\bibitem{Rotar2} V. I. Rotar' (1979). Limit theorems for polylinear forms. \textit{J. Multivariate Anal.} {\bf 9}, 511--530.

\bibitem{Schatte}
P. Schatte (1988).
On strong versions of the central limit theorem.
{\it Math. Nachr.} {\bf 137}, 249-256.

\bibitem{Stanley} R. Stanley (1997). \textit{Enumerative combinatorics, Vol.
1. }Cambridge University Press.


\bibitem{TV} T. Tao and V. Vu (2008).
Random matrices: Universality of ESD and the Circular Law
(with an appendix by M. Krishnapur).
To appear in: {\it Ann. Probab.}


\end{thebibliography}
\end{document}